\RequirePackage{txfonts}\normalfont     

\documentclass[11pt,%
hypertex,%
onecolumn,oneside,a4paper,article
,italian,frenchb,german,
swedish,british%
]{memoir}

\newcommand{\usetxfonts}{\usepackage{txfonts}}

\newcommand*\savesymbol[1]{%
  \expandafter\let\csname orig#1\expandafter\endcsname\csname#1\endcsname
  \expandafter\let\csname #1\endcsname\relax
}

\newcommand*\restoresymbol[2]{%
  \expandafter\global\expandafter\let\csname#1#2\expandafter\endcsname%
    \csname#2\endcsname
  \expandafter\global\expandafter\let\csname#2\expandafter\endcsname%
    \csname orig#2\endcsname
}

\usepackage[T1]{fontenc} 
\usepackage{textcomp}
\usepackage{amsmath}
\usepackage{amssymb}
\usepackage{amsxtra}
\usepackage[
]{babel}
\newcommand{\langfrench}{\foreignlanguage{french}}

\newcommand{\langitalian}{\foreignlanguage{italian}}

\newcommand{\langlatin}{\foreignlanguage{latin}}

\savesymbol{qedsymbol}

\usepackage{amsthm}

\restoresymbol{ams}{qedsymbol}
\newcommand{\QEE}
{\amsqedsymbol}
\newcommand{\qee}{\leavevmode\unskip\penalty9999 \hbox{}\nobreak\hfill
\quad\hbox{\QEE}}

\providecommand{\firstname}[1]{#1}
\providecommand{\surname}[1]{#1}

\providecommand{\usetxfonts}{}
\usetxfonts

\newcommand{\Alpha}{\mathrm{A}}
\let\varAlpha A
\newcommand{\Beta}{\mathrm{B}}
\let\varBeta B
\newcommand{\Epsilon}{\mathrm{E}}
\let\varEpsilon E
\newcommand{\Zeta}{\mathrm{Z}}
\let\varZeta Z
\newcommand{\Eta}{\mathrm{H}}
\let\varEta H
\newcommand{\Iota}{\mathrm{I}}
\let\varIota I
\newcommand{\Kappa}{\mathrm{K}}
\let\varKappa K
\newcommand{\Mu}{\mathrm{M}}
\let\varMu M
\newcommand{\Nu}{\mathrm{N}}
\let\varNu N
\newcommand{\Omicron}{\mathrm{O}}
\let\varOmicron O
\newcommand{\Rho}{\mathrm{P}}
\let\varRho P
\newcommand{\Tau}{\mathrm{T}}
\let\varTau T
\newcommand{\Chi}{\mathrm{X}}
\let\varChi X
\newcommand{\omicron}{\mathrm{o}}

\DeclareFontFamily{U}{egreek}{\skewchar\font'177}%
\DeclareFontShape{U}{egreek}{m}{n}{<-6>s*[1]eurm5<6-8>s*[1]eurm7 <8->s*[1]eurm10}{}%
\DeclareFontShape{U}{egreek}{m}{it}{<->s*[1]eurmo10}{}%
\DeclareFontShape{U}{egreek}{b}{n}{<-6>s*[1]eurb5 <6-8>s*[1]eurb7 <8->s*[1]eurb10}{}%
\DeclareFontShape{U}{egreek}{b}{it}{<->s*[1]eurbo10}{}%
\DeclareSymbolFont{egreeki}{U}{egreek}{m}{it}%
\SetSymbolFont{egreeki}{bold}{U}{egreek}{b}{it}
\DeclareMathSymbol{\epartial}{\mathalpha}{egreeki}{"40}

\DeclareMathSymbol{\ealpha}{\mathalpha}{egreeki}{"0B}
\DeclareMathSymbol{\ebeta}{\mathalpha}{egreeki}{"0C}
\DeclareMathSymbol{\egamma}{\mathalpha}{egreeki}{"0D}
\DeclareMathSymbol{\edelta}{\mathalpha}{egreeki}{"0E}
\DeclareMathSymbol{\eepsilon}{\mathalpha}{egreeki}{"0F}
\DeclareMathSymbol{\ezeta}{\mathalpha}{egreeki}{"10}
\DeclareMathSymbol{\eeta}{\mathalpha}{egreeki}{"11}
\DeclareMathSymbol{\etheta}{\mathalpha}{egreeki}{"12}
\DeclareMathSymbol{\eiota}{\mathalpha}{egreeki}{"13}
\DeclareMathSymbol{\ekappa}{\mathalpha}{egreeki}{"14}
\DeclareMathSymbol{\elambda}{\mathalpha}{egreeki}{"15}
\DeclareMathSymbol{\emu}{\mathalpha}{egreeki}{"16}
\DeclareMathSymbol{\enu}{\mathalpha}{egreeki}{"17}
\DeclareMathSymbol{\exi}{\mathalpha}{egreeki}{"18}
\DeclareMathSymbol{\eomicron}{\mathalpha}{egreeki}{"6F}
\DeclareMathSymbol{\epi}{\mathalpha}{egreeki}{"19}
\DeclareMathSymbol{\erho}{\mathalpha}{egreeki}{"1A}
\DeclareMathSymbol{\esigma}{\mathalpha}{egreeki}{"1B}
\DeclareMathSymbol{\etau}{\mathalpha}{egreeki}{"1C}
\DeclareMathSymbol{\eupsilon}{\mathalpha}{egreeki}{"1D}
\DeclareMathSymbol{\ephi}{\mathalpha}{egreeki}{"1E}
\DeclareMathSymbol{\echi}{\mathalpha}{egreeki}{"1F}
\DeclareMathSymbol{\epsi}{\mathalpha}{egreeki}{"20}
\DeclareMathSymbol{\eomega}{\mathalpha}{egreeki}{"21}
\DeclareMathSymbol{\evarepsilon}{\mathalpha}{egreeki}{"22}
\DeclareMathSymbol{\evartheta}{\mathalpha}{egreeki}{"23}
\DeclareMathSymbol{\evarpi}{\mathalpha}{egreeki}{"24}
\let\evarrho\erho 
\let\evarsigma\esigma
\let\evarkappa\ekappa
\DeclareMathSymbol{\evarphi}{\mathalpha}{egreeki}{"27}
\DeclareMathSymbol{\evarAlpha}{\mathalpha}{egreeki}{"41}
\DeclareMathSymbol{\evarBeta}{\mathalpha}{egreeki}{"42}
\DeclareMathSymbol{\evarGamma}{\mathalpha}{egreeki}{"00}
\DeclareMathSymbol{\evarDelta}{\mathalpha}{egreeki}{"01}
\DeclareMathSymbol{\evarEpsilon}{\mathalpha}{egreeki}{"45}
\DeclareMathSymbol{\evarZeta}{\mathalpha}{egreeki}{"5A}
\DeclareMathSymbol{\evarEta}{\mathalpha}{egreeki}{"48}
\DeclareMathSymbol{\evarTheta}{\mathalpha}{egreeki}{"02}
\DeclareMathSymbol{\evarIota}{\mathalpha}{egreeki}{"49}
\DeclareMathSymbol{\evarKappa}{\mathalpha}{egreeki}{"4B}
\DeclareMathSymbol{\evarLambda}{\mathalpha}{egreeki}{"03}
\DeclareMathSymbol{\evarMu}{\mathalpha}{egreeki}{"4D}
\DeclareMathSymbol{\evarNu}{\mathalpha}{egreeki}{"4E}
\DeclareMathSymbol{\evarXi}{\mathalpha}{egreeki}{"04}
\DeclareMathSymbol{\evarOmicron}{\mathalpha}{egreeki}{"4F}
\DeclareMathSymbol{\evarPi}{\mathalpha}{egreeki}{"05}
\DeclareMathSymbol{\evarRho}{\mathalpha}{egreeki}{"50}
\DeclareMathSymbol{\evarSigma}{\mathalpha}{egreeki}{"06}
\DeclareMathSymbol{\evarTau}{\mathalpha}{egreeki}{"54}
\DeclareMathSymbol{\evarUpsilon}{\mathalpha}{egreeki}{"07}
\DeclareMathSymbol{\evarPhi}{\mathalpha}{egreeki}{"08}
\DeclareMathSymbol{\evarChi}{\mathalpha}{egreeki}{"58}
\DeclareMathSymbol{\evarPsi}{\mathalpha}{egreeki}{"09}
\DeclareMathSymbol{\evarOmega}{\mathalpha}{egreeki}{"0A} 
\DeclareSymbolFont{egreekr}{U}{egreek}{m}{n}%
\SetSymbolFont{egreekr}{bold}{U}{egreek}{b}{n}
\DeclareMathSymbol{\epartialup}{\mathalpha}{egreekr}{"40}

\DeclareMathSymbol{\ealphaup}{\mathalpha}{egreekr}{"0B}
\DeclareMathSymbol{\ebetaup}{\mathalpha}{egreekr}{"0C}
\DeclareMathSymbol{\egammaup}{\mathalpha}{egreekr}{"0D}
\DeclareMathSymbol{\edeltaup}{\mathalpha}{egreekr}{"0E}
\DeclareMathSymbol{\eepsilonup}{\mathalpha}{egreekr}{"0F}
\DeclareMathSymbol{\ezetaup}{\mathalpha}{egreekr}{"10}
\DeclareMathSymbol{\eetaup}{\mathalpha}{egreekr}{"11}
\DeclareMathSymbol{\ethetaup}{\mathalpha}{egreekr}{"12}
\DeclareMathSymbol{\eiotaup}{\mathalpha}{egreekr}{"13}
\DeclareMathSymbol{\ekappaup}{\mathalpha}{egreekr}{"14}
\DeclareMathSymbol{\elambdaup}{\mathalpha}{egreekr}{"15}
\DeclareMathSymbol{\emuup}{\mathalpha}{egreekr}{"16}
\DeclareMathSymbol{\enuup}{\mathalpha}{egreekr}{"17}
\DeclareMathSymbol{\exiup}{\mathalpha}{egreekr}{"18}
\DeclareMathSymbol{\eomicronup}{\mathalpha}{egreekr}{"6F}
\DeclareMathSymbol{\epiup}{\mathalpha}{egreekr}{"19}
\DeclareMathSymbol{\erhoup}{\mathalpha}{egreekr}{"1A}
\DeclareMathSymbol{\esigmaup}{\mathalpha}{egreekr}{"1B}
\DeclareMathSymbol{\etauup}{\mathalpha}{egreekr}{"1C}
\DeclareMathSymbol{\eupsilonup}{\mathalpha}{egreekr}{"1D}
\DeclareMathSymbol{\ephiup}{\mathalpha}{egreekr}{"1E}
\DeclareMathSymbol{\echiup}{\mathalpha}{egreekr}{"1F}
\DeclareMathSymbol{\epsiup}{\mathalpha}{egreekr}{"20}
\DeclareMathSymbol{\eomegaup}{\mathalpha}{egreekr}{"21}
\DeclareMathSymbol{\evarepsilonup}{\mathalpha}{egreekr}{"22}
\DeclareMathSymbol{\evarthetaup}{\mathalpha}{egreekr}{"23}
\DeclareMathSymbol{\evarpiup}{\mathalpha}{egreekr}{"24}
\let\evarrhoup\erhoup 
\let\evarsigmaup\esigmaup
\let\evarkappaup\ekappaup
\DeclareMathSymbol{\evarphiup}{\mathalpha}{egreekr}{"27}
\DeclareMathSymbol{\eAlpha}{\mathalpha}{egreekr}{"41}
\DeclareMathSymbol{\eBeta}{\mathalpha}{egreekr}{"42}
\DeclareMathSymbol{\eGamma}{\mathalpha}{egreekr}{"00}
\DeclareMathSymbol{\eDelta}{\mathalpha}{egreekr}{"01}
\DeclareMathSymbol{\eEpsilon}{\mathalpha}{egreekr}{"45}
\DeclareMathSymbol{\eZeta}{\mathalpha}{egreekr}{"5A}
\DeclareMathSymbol{\eEta}{\mathalpha}{egreekr}{"48}
\DeclareMathSymbol{\eTheta}{\mathalpha}{egreekr}{"02}
\DeclareMathSymbol{\eIota}{\mathalpha}{egreekr}{"49}
\DeclareMathSymbol{\eKappa}{\mathalpha}{egreekr}{"4B}
\DeclareMathSymbol{\eLambda}{\mathalpha}{egreekr}{"03}
\DeclareMathSymbol{\eMu}{\mathalpha}{egreekr}{"4D}
\DeclareMathSymbol{\eNu}{\mathalpha}{egreekr}{"4E}
\DeclareMathSymbol{\eXi}{\mathalpha}{egreekr}{"04}
\DeclareMathSymbol{\eOmicron}{\mathalpha}{egreekr}{"4F}
\DeclareMathSymbol{\ePi}{\mathalpha}{egreekr}{"05}
\DeclareMathSymbol{\eRho}{\mathalpha}{egreekr}{"50}
\DeclareMathSymbol{\eSigma}{\mathalpha}{egreekr}{"06}
\DeclareMathSymbol{\eTau}{\mathalpha}{egreekr}{"54}
\DeclareMathSymbol{\eUpsilon}{\mathalpha}{egreekr}{"07}
\DeclareMathSymbol{\ePhi}{\mathalpha}{egreekr}{"08}
\DeclareMathSymbol{\eChi}{\mathalpha}{egreekr}{"58}
\DeclareMathSymbol{\ePsi}{\mathalpha}{egreekr}{"09}
\DeclareMathSymbol{\eOmega}{\mathalpha}{egreekr}{"0A}

\newcommand{\alleugreek}{%
\let\alpha\ealpha
\let\beta\ebeta
\let\gamma\egamma
\let\delta\edelta
\let\epsilon\eepsilon
\let\zeta\ezeta
\let\eta\eeta
\let\theta\etheta
\let\iota\eiota
\let\kappa\ekappa
\let\lambda\elambda
\let\mu\emu
\let\nu\enu
\let\xi\exi
\let\omicron\eomicron
\let\pi\epi
\let\rho\erho
\let\sigma\esigma
\let\tau\etau
\let\upsilon\eupsilon
\let\phi\ephi
\let\chi\echi
\let\psi\epsi
\let\omega\eomega
\let\varepsilon\evarepsilon
\let\vartheta\evartheta
\let\varpi\evarpi
\let\varrho\evarrho 
\let\varsigma\evarsigma
\let\varkappa\evarkappa
\let\varphi\evarphi
\let\varAlpha\evarAlpha
\let\varBeta\evarBeta
\let\varGamma\evarGamma
\let\varDelta\evarDelta
\let\varEpsilon\evarEpsilon
\let\varZeta\evarZeta
\let\varEta\evarEta
\let\varTheta\evarTheta
\let\varIota\evarIota
\let\varKappa\evarKappa
\let\varLambda\evarLambda
\let\varMu\evarMu
\let\varNu\evarNu
\let\varXi\evarXi
\let\varOmicron\evarOmicron
\let\varPi\evarPi
\let\varRho\evarRho
\let\varSigma\evarSigma
\let\varTau\evarTau
\let\varUpsilon\evarUpsilon
\let\varPhi\evarPhi
\let\varChi\evarChi
\let\varPsi\evarPsi
\let\varOmega\evarOmega 
\let\alphaup\ealphaup
\let\betaup\ebetaup
\let\gammaup\egammaup
\let\deltaup\edeltaup
\let\epsilonup\eepsilonup
\let\zetaup\ezetaup
\let\etaup\eetaup
\let\thetaup\ethetaup
\let\iotaup\eiotaup
\let\kappaup\ekappaup
\let\lambdaup\elambdaup
\let\muup\emuup
\let\nuup\enuup
\let\xiup\exiup
\let\omicronup\eomicronup
\let\piup\epiup
\let\rhoup\erhoup
\let\sigmaup\esigmaup
\let\tauup\etauup
\let\upsilonup\eupsilonup
\let\phiup\ephiup
\let\chiup\echiup
\let\psiup\epsiup
\let\omegaup\eomegaup
\let\varepsilonup\evarepsilonup
\let\varthetaup\evarthetaup
\let\varpiup\evarpiup
\let\varrhoup\evarrhoup 
\let\varsigmaup\evarsigmaup
\let\varkappaup\evarkappaup
\let\varphiup\evarphiup
\let\Alpha\eAlpha
\let\Beta\eBeta
\let\Gamma\eGamma
\let\Delta\eDelta
\let\Epsilon\eEpsilon
\let\Zeta\eZeta
\let\Eta\eEta
\let\Theta\eTheta
\let\Iota\eIota
\let\Kappa\eKappa
\let\Lambda\eLambda
\let\Mu\eMu
\let\Nu\eNu
\let\Xi\eXi
\let\Omicron\eOmicron
\let\Pi\ePi
\let\Rho\eRho
\let\Sigma\eSigma
\let\Tau\eTau
\let\Upsilon\eUpsilon
\let\Phi\ePhi
\let\Chi\eChi
\let\Psi\ePsi
\let\Omega\eOmega
}

\newcommand{\unambeugreek}{%
\let\omicron\eomicron
\let\varAlpha\evarAlpha
\let\varBeta\evarBeta
\let\varEpsilon\evarEpsilon
\let\varZeta\evarZeta
\let\varEta\evarEta
\let\varIota\evarIota
\let\varKappa\evarKappa
\let\varMu\evarMu
\let\varNu\evarNu
\let\varOmicron\evarOmicron
\let\varRho\evarRho
\let\varTau\evarTau
\let\varChi\evarChi
\let\omicronup\eomicronup
\let\Alpha\eAlpha
\let\Beta\eBeta
\let\Epsilon\eEpsilon
\let\Zeta\eZeta
\let\Eta\eEta
\let\Iota\eIota
\let\Kappa\eKappa
\let\Mu\eMu
\let\Nu\eNu
\let\Omicron\eOmicron
\let\Rho\eRho
\let\Tau\eTau
\let\Chi\eChi
}

\usepackage{bm}
\providecommand{\piup}{\epiup}
\providecommand{\deltaup}{\edeltaup}
\providecommand{\zetaup}{\ezetaup}
\providecommand{\varepsilonup}{\evarepsilonup}

\providecommand{\coloneqq}{\mathrel{\mathop:}=}

\newcommand{\delt}{\deltaup}
\newcommand{\I}{\mathrm{i}}
\newcommand{\di}{\mathrm{d}}

\DeclareSymbolFont{AMSb}{U}{msb}{m}{n}
\DeclareMathSymbol{\CC}{\mathalpha}{AMSb}{"43}
\DeclareMathSymbol{\RR}{\mathalpha}{AMSb}{"52}
\DeclareMathSymbol{\NN}{\mathalpha}{AMSb}{"4E}
\DeclareMathSymbol{\ZZ}{\mathalpha}{AMSb}{"5A}
\DeclareMathSymbol{\QQ}{\mathalpha}{AMSb}{"51}

\expandafter\let\csname im\expandafter\endcsname\csname Im\endcsname
  \expandafter\let\csname Im\endcsname\relax

\expandafter\let\csname re\expandafter\endcsname\csname Re\endcsname
  \expandafter\let\csname Re\endcsname\relax
\DeclareMathOperator{\Re}{Re}

\newcommand{\defd}{\coloneqq}

\newcommand{\Lor}{\bigvee}
\newcommand{\lonlyif}{\mathbin{\Rightarrow}}

\newcommand{\cond}
{\mathpunct{|}}
\newcommand{\st}{\mid}
\newcommand{\inn}{\cdot}
\newcommand{\dotv}{\mathord{\,\cdot\,}}

\renewcommand{\le}{\leqslant}
\renewcommand{\ge}{\geqslant}
\DeclareMathDelimiter{\lclose}{\mathopen}{operators}{"5B}{largesymbols}{"02}
\DeclareMathDelimiter{\rclose}{\mathclose}{operators}{"5D}{largesymbols}{"03}
\DeclareMathDelimiter{\lopen}{\mathopen}{operators}{"5D}{largesymbols}{"03}
\DeclareMathDelimiter{\ropen}{\mathclose}{operators}{"5B}{largesymbols}{"02}
\newcommand{\clcl}[1]{\lclose#1\rclose}

\newcommand{\abs}[1]{\lvert#1\rvert}

\newcommand{\set}[1]{\{#1\}}
\DeclareMathOperator{\pr}{P}
\DeclareMathOperator{\G}{\Gamma}
\newcommand{\sect}{\S~}
\newcommand{\sects}{\S\S~}
\newcommand{\chap}{ch.~}%
\newcommand{\chaps}{chs~}%
\newcommand{\eqn}{eq.~}%
\theoremstyle{remark}
\newtheorem{remark}{Remark}
\theoremstyle{definition}

\newcommand{\ie}{{i.e.}}

\newcommand{\eg}{{e.g.}}

\newcommand{\cf}{{cf.}}
\newcommand{\Cf}{{Cf.}}

\newcommand{\etal}{{et al.}}

\newcommand{\bd}{\hspace{0pt}}%

\newcommand{\labelbis}[1]{\tag*{(\ref{#1})$_\text{r}$}}

\providecommand{\href}[2]{#2}

\DeclareMathOperator{\p}{P}
\renewcommand{\|}{\cond}

\expandafter\let\csname amp\expandafter\endcsname\csname&\endcsname
  \expandafter\let\csname &\endcsname\relax

  \expandafter\global\expandafter\let\csname&\expandafter\endcsname%
    \csname land\endcsname

\alleugreek

\usepackage{paralist}

\renewcommand{\langlatin}{\foreignlanguage{nohyphenation}}

\usepackage[citestyle=numeric-comp,%
bibstyle=numericmana,%
sorting=none,%
natbib=false,firstinits=true,%
block=space,%
hyperref=true,%
useprefix=true
]{biblatex}
\renewcommand{\multicitedelim}{\addsemicolon\space}

\bibliography{mana090501-PT_MEv2.bbl}

\renewcommand{\bibfont}{\small}

\usepackage[misc]{ifsym}

\chapterstyle{reparticle}
\checkandfixthelayout
\fixpdflayout

\usepackage[%
dvips,
breaklinks,pdfauthor={P G L  Porta Mana}]{hyperref}

\hypersetup{pdfauthor={P G L  Porta Mana},pdftitle={Plausibility logic and
    maximum entropy}}

\usepackage{graphicx}

\setlxvchars
\setxlvchars

\providecommand{\affiliation}[1]{\textit{\small #1}}
\providecommand{\pacs}[1]{{\footnotesize\textsc{PACS} numbers: #1}}
\providecommand{\msc}[1]{{\footnotesize\textsc{MSC} numbers: #1}}
\providecommand{\email}[1]{\texttt{\href{mailto:#1}{#1}}}

\newcommand{\asudedication}[1]{%
{\par\centering\textit{#1}\par}}

\title{On the relation between plausibility logic and the maximum-entropy
  principle: a numerical study}

\author{\firstname{P. G. L.}\;\surname{Porta Mana},\thanks{Email:
    \email{lmana@zperimeterinstitute.ca} (remove the z)}
\\
\affiliation{Perimeter Institute for Theoretical Physics, Canada}
}

\date{10 November 2009}

\copypagestyle{manaart}{plain}
\makeheadrule{manaart}{\headwidth}{\normalrulethickness}
\makeoddhead{manaart}{\small\textsc{Porta
    Mana}}{}{\small\textit{Plausibility logic and maximum entropy}}

\newenvironment{acknowledgements}{\chapter*{Acknowledgements}\addcontentsline{toc}{chapter}{Acknowledgements}}{\par}

\theoremstyle{plain}

\newcommand{\tprod}{\mathop{\textstyle\prod}\limits}
\newcommand{\tsum}{\mathop{\textstyle\sum}\limits}

\newcommand{\lsum}{\mathop{\textstyle\sum}\nolimits}
\newcommand{\lprod}{\mathop{\textstyle\prod}\nolimits}

\newcommand{\me}{maxi\-mum-\bd entropy}

\newcommand{\yG}{g}

\newcommand{\yU}{\text{\scriptsize\Cube{1}}}
\newcommand{\yB}{\text{\scriptsize\Cube{2}}}
\newcommand{\yT}{\text{\scriptsize\Cube{3}}}
\newcommand{\yQ}{\text{\scriptsize\Cube{4}}}
\newcommand{\yC}{\text{\scriptsize\Cube{5}}}
\newcommand{\yS}{\text{\scriptsize\Cube{6}}}

\newcommand{\yk}{K}
\newcommand{\ykk}{m}
\newcommand{\ykkk}{\bm{\ykk}}
\newcommand{\yl}{L}
\newcommand{\yD}{D}
\newcommand{\ybu}{H_{\text{B}}}

\newcommand{\yA}{A}
\newcommand{\yf}{f}
\newcommand{\yff}{\bm{\yf}}
\newcommand{\ydf}{\di\bm{\yf}}
\newcommand{\yp}{\bm{p}}
\newcommand{\ydp}{\di\yp}
\newcommand{\iid}{i.i.d.}

\newcommand{\yiid}{I_{\text{ft}}}
\newcommand{\ycom}[1]{I_{\text{m}}^{#1}}
\newcommand{\ydir}[1]{I_{\text{J}}^{#1}}
\newcommand{\yn}{\bm{N}}
\newcommand{\ync}{c}
\newcommand{\yden}{n!\, \ydp}
\newcommand{\ydens}{\bm{\omega}}
\newcommand{\ysimp}{\Delta}
\newcommand{\ysimpa}{{\Delta_a}}

\newcommand{\we}{\wedge}
\newcommand{\yo}{\bm{v}}
\newcommand{\yR}[2]{R^{(#1)}_{#2}}
\newcommand{\asy}{\simeq}
\newcommand{\ycm}{\varepsilon}

\DeclareBibliographyCategory{extra}
\DeclareBibliographyCategory{extrb}
\DeclareBibliographyCategory{extrc}
\DeclareBibliographyCategory{extrd}
\DeclareBibliographyCategory{extre}
\DeclareBibliographyCategory{6}

\newcommand{\citep}{\parencites}
\newcommand{\citey}{\parencites*}
\newcommand{\citein}[2][]{\textnormal{\textcite[#1]{#2}}\addtocategory{extra}{#2}}
\newcommand{\citebi}[1]{\textcite{#1}\addtocategory{extra}{#1}}
\renewcommand{\cite}{\citep}

\begin{document}
\setlength{\droptitle}{-3\onelineskip}




\maketitle

\abslabeldelim{:\quad}
\setlength{\abstitleskip}{-\absparindent}
\abstractrunin
\begin{abstract}
  What is the relationship between plausibility logic and the principle of
  maximum entropy? When does the principle give unreasonable or wrong
  results? When is it appropriate to use the rule `$\text{expectation} =
  \text{average}$'? Can plausibility logic give the same answers as the
  principle, and better answers if those of the principle are unreasonable?
  To try to answer these questions, this study offers a numerical
  collection of plausibility distributions given by the \me\ principle and
  by plausibility logic for a set of fifteen simple problems: throwing dice.
  \\[2\jot]
  \pacs{02.50.Cw,02.50.Tt,01.70.+w}\\
  \msc{03B48,60G09,60A05}
\end{abstract}
\asudedication{\langitalian{Dedicato a mia madre per il suo trentesimo compleanno}}


\begin{refsegment}
\chapter{When and how should the \me\ principle be applied?}
\label{cha:maxent_plaus}

For the student of plausibility logic\footnote{I call `plausibility logic'
  what many other authors call `(Bayesian) probability theory'. `Logic',
  because it is a generalization of the truth-logical calculus.
  `Plausibility', because `degree of belief' is unfortunately too unwieldy
  and many authors still contend that $\text{`probability'} =
  \text{`frequency'}$ or, perhaps worse, $\text{`probability'} =
  \text{`(Lebesgue) measure'}$
.}, the theory of the principles governing plausible inference, the
application of the theory in any given problem is crystal clear in
principle:
\begin{inparaenum}[(1)]
\item 
The problem is analysed and reduced to a set of
propositions $\set{A_i}$ and background knowledge $I$.
\item Some plausibilities $\p[A_{i_1} \dotso A_{i_2} \| (A_{i_3} \dotso
  A_{i_4}) \&I] \in \clcl{0,1}$ are assigned, consistently with the laws
  below, according to our actual or hypothetical knowledge of the situation
  and to convenience; the `$A_{i_j} \dotso A_{i_k}$' represent collections
  of $A_i$s joined by various logical connectives (`$\lnot$', `$\land$',
  `$\lor$', `$\lonlyif$').
\item Finally, using the basic laws
  \begin{subequations}
    \begin{align}
      \label{eq:basic_laws}
      \p(\lnot A_i \| I) &= 1-\p(A_i \| I),
      \\
      \p(A_i \land A_j \| I) &= \p(A_i \| A_j \&I)\, \p(A_j \|I),
      \\
      \p(A_i \lor A_j \| I) &= \p(A_i \|I) + \p(A_j \|I) - \p(A_i \land A_j
      \| I),
      \\
      \p(A_i \lonlyif A_j \| I) &= \p(\lnot A_i \|I) +
\p(A_j \| A_i \&I)\, \p(A_i \|I),
    \end{align}
  \end{subequations}
  the values of unassigned plausibilities of interest, or at least some
  bounds for them, are calculated --- \eg\ via fractional linear programming,
  as explained in the scholarly and unfortunately much neglected works
  of Hailperin \citey{hailperin1984}[\sects0.4, 4.5, 5.4, 6.2 and
  passim]{hailperin1996}{hailperin2006}. If required, some mathematical
  limits are taken.
\end{inparaenum}
This procedure applies whether we want to calculate the plausibility of a
proposition of interest, to explore the plausible consequences of a
hypothesis, or to see the relationship between the plausibilities assigned
upon different contexts; and other similar problems. 

When the same student approaches the principle of maximum entropy, he is
not welcomed by a comparable level of clarity. For example, he is told to
use data consisting in an average value; but does the number of
observations contributing to this average matter? He is told to equate an
average with an expectation; but when and why is it reasonable to equate
these very different quantities?
He is told that the principle yields a plausibility distribution; but does
this plausibility concern one of the observations contributing to the data,
or a new observation? And what is the conditional of the probability
distribution given by the principle?\ that is, if the principle gives
a $\pr(A_i \cond \dotv)$, what proposition does the dot stand
for? In fact, the specification of the conditional, or `supposal'
\citep{johnson1924,johnson1932,johnson1932b,johnson1932c} or `context', of
a plausibility is extremely important, also because the conditionals and
the arguments of two plausibilities must match in a precise way for the
plausibility laws to be used. It is meaningless to multiply, \eg, $\p(A_3
\|I)$ and $\p(A_2 \|A_1 \&I)$ without further qualifications.

The various `proofs' of the maximum-entropy principle in the literature do
not help the student very much either. The most rigorous of them cover only
specific situations, leaving out important ones. For example,
van Campenhout \amp\ Cover \citey{vancampenhoutetal1981} and Csisz\'ar
\citey{csiszar1985} prove that a uniform \iid\ model gives, to any
observation contributing to the average, the same distribution as the \me\
principle, when the average used as data comes from a very
large number of observations. But they are silent on whether that
distribution is valid for similar observations \emph{not} contributing to
the average. On the other hand, the logics behind many proofs covering this
last application of the principle have been repeatedly attacked by various
authors. I recommend Skyrms' \citey{skyrms1987} and especially Uffink's
\citey{uffink1996b} analyses, which give insights on several important points
(but show confusion about others), many of which are repeated here.

One could even say that there is not just one `\me\ principle', because in
the literature this term denotes \emph{qualitatively} different procedures
and problems, in which one seeks qualitatively different kinds of
distributions --- \eg, distributions of probability, of frequency, or of
intensity as in image-reconstruction problems. In some of these cases there
really is no `principle' but only a `rule' which appears asymptotically
from choosing the frequency- or intensity-distribution
$M_i/M$, $i=1,\dotsc,r$, that can be realized in most ways.
The number of ways is usually given by the multiplicity factor
\begin{equation}
  \label{eq:multfactor}
  \frac{M!}{\prod_{i=1}^r M_i!} 
=  \ycm(M_i)\, \exp[M\, H(M_i/M)],
\quad \text{with $(M+1)^{-r} \le \ycm(M_i) \le 1$},
\end{equation}
and this is asymptotically equal
\citep[\sect1.2]{csiszaretal1981}[\sect2.1]{csiszaretal2004b} to the
exponential of a very large multiple of the Shannon entropy of the
distribution,
\begin{equation}
  \label{eq:shan-entropy}
  H(\yf_i) \defd -\lsum_i \yf_i \ln \yf_i,
\qquad\text{$\yf_i \ge 0$, $\lsum_i \yf_i=1$; $0 \ln 0 \defd 0$}.
\end{equation}
The distribution chosen by `counting reasons' is therefore the one having
maximum Shannon entropy. In examples like these the \me\ rule is an obvious
consequence of plausibility logic with an assumption of symmetry on a
particular hypothesis space (\eg, the set of outcome sequences) and does
not need additional principles for its justification. These cases do not
concern us here. But the counting arguments above make no sense for
distributions of plausibility, and the application of the rule, especially
to \emph{new} observations, seems to require additional principles besides
those of plausibility logic.



\chapter{A numerical comparison: fifteen problems}
\label{cha:comparison}

The purpose of the present work is to examine the \me\ principle `in
action' in a collection of simple problems, to see under which circumstances
its results are intuitively reasonable or not and how these compare with
those given by plausibility logic alone.

The collection of problems is the following: Of $N$ throws of a die we know
the average $a$ and nothing else; not even the throwing technique or the
kind of die used, although we assume them to be the same or at least very
similar in all throws. We want the plausibility distributions for:
\begin{enumerate}
\item the outcomes of one of the $N$ throws contributing to the average,
  which we call `old throws'; and
\item the outcomes of a throw --- of the same or very similar die and with
  the same throwing technique --- outside the set of $N$ old throws, which
  we call a `new throw'.
\end{enumerate}
We consider the problems obtained by combining the particular values $N=1$,
$6$, $12$ and $N$ large, together with $a=6$, $5$, and $7/2$,
for a total of fifteen problems.

The adjectives `old' and `new' qualifying the throws have really no temporal
meaning; in fact, the problem is \emph{by assumption} completely symmetric
with respect to the temporal ordering of all throws, or the way they are
scattered in space-time. A `new' throw could precede an `old' one, or they
could all happen at the same time. To stress this I shall use the present
tense even with `old' throws, saying \eg\ `the second old throw
\emph{gives} face \yS'. Plausibility logic is affected by \emph{what} we
know; it is immaterial \emph{how} and \emph{when} we come to know it (even
if it had been by precognition or other imaginary ways).

The old throws are numbered from 1 to $N$ (this numbering, again, bearing
no temporal meaning); a new throw is assigned the number $0$. The
proposition stating that the outcome of throw $j$ is the face `$i$' of the
die is denoted by $\yR{j}{i}$. The proposition stating that the average of
the $N$ throws is $a$ is denoted by $\yA^N_a$. Thus the plausibility
distributions sought are, in symbols,
\begin{enumerate}
\item $\p(\yR{1}{i} \| \yA^N_a \land I)$\quad(old throw, $j=1$) and
\item $\p(\yR{0}{i} \| \yA^N_a \land I)$\quad(new throw, $j=0$).
\end{enumerate}

The proposition $I$ states our background knowledge, \ie, mathematically
speaking, our plausibility model. An important question indeed is: what kind
of background knowledge, or plausibility model, are we implicitly assuming
when we use the \me\ principle?

\section{Exchangeable plausibility models used for comparison}
\label{cha:models}

In our comparison we consider three plausibility models, all three
infinitely exchangeable. Infinite exchangeability means that, in assigning
a plausibility distribution to the set of outcomes of \emph{any} number of
old and new throws, we do not care which throw each outcome belongs to. An alternative
interpretation is the following: if we knew enough about the circumstances
of each throw, knowledge of all other throws would be irrelevant in our
plausibility assignment for that throw; but those circumstances are unknown
and we can only assign a plausibility distribution to their various
possibilities, which are themselves plausibility-indexed; see
\citep{johnson1932c,caves2000c,portamanaetal2007}. In both interpretations
the plausibility distribution for the outcomes of throws $j_1, \dotsc, j_n$
has the form
\begin{equation}
  \label{eq:exchangeable}
  \p(R^{(j_1)}_{i_n} \& \dotsb \& R^{(j_1)}_{i_n} \| I)
=
\int p_{i_1} \dotsm p_{i_n} \, \yG(\yp\| I)\,\ydp;
\end{equation}
the integration is over the simplex
of plausibility distributions
\begin{equation}
  \label{eq:simplex}
  \ysimp \defd 
\set{\yp \defd (p_1, \dotsc, p_6) 
\st p_i \ge 0, \lsum_i p_i = 1}.
\end{equation}
The generalized function \citep{egorov1990,egorov1990b,egorov2001} $\yG$
characterizes the exchangeable model chosen; $\yG\,\ydp$ can be interpreted
as the plausibility density for the limiting frequencies of the outcomes as
the number of new throws increases indefinitely or, in the alternative
interpretation of infinite exchangeability, as the plausibility density for
those circumstances with index values around the volume element $\ydp$. See
app.~\ref{cha:can_dens} for some remarks about this volume element.

The choice of $\yG(\yp \| I)$ defines our exchangeable model $I$. For
example, for a generalized density concentrated on the uniform
distribution,
\begin{equation}
  \label{eq:G-iid}
 \yG(\yp \| \yiid)\,\ydp \defd \tprod_i \delt(p_i - 1/6)\,\ydp,
\end{equation}
the model $\yiid$ gives each throw a uniform plausibility distribution,
independent of the knowledge about all other throws and identical for all
throws, `\iid'. It represents the unshakeable belief that the die and
throwing technique be absolutely `fair'.
We call this the \emph{fair-throw} model; it will be the first to be
compared with the \me\ principle.

Suppose, on the other hand, that we judge the plausibility for face `$i$'
on a new throw to depend only on $N$ and on how many old throws yield face
`$i$'; in other words, the frequencies of the other faces of old throws are
irrelevant. Then we have the \emph{Johnson} model, also called
\emph{Dirichlet} model, $\ydir{\yk}$. It depends on a parameter $\yk$ and is
represented by the density
\begin{equation}
  \label{eq:G-dir}
  \yG(\yp \| \ydir{\yk})\,\ydp \defd 
\G(1+\lsum_i\yk)\,\Biggl[\prod_i \frac{p_i^{\yk-1}}{\G(\yk)}\Biggr]\ydp,
\qquad \yk >0;
\end{equation}
see Johnson \citey{johnson1924,johnson1932c}, Good
\citey[\chap4]{good1965}, Zabell \citey{zabell1982}, Jaynes
\citey{jaynes1986d_r1996}, and references therein. For the values $\yk=1,
2, 5, 50$ and $\yk$ large this will be the second model in our numerical
comparison. The case $\yk=1$ (a constant density) corresponds to the
multidimensional form of Bayes' \citey[Scholium]{bayes1763} and Laplace's
suggestion \citep[p.~\textsc{xvii}]{laplace1814_r1819}: `\langfrench{Quand
  la probabilit\'e d'un \'ev\'enement simple est inconnue, on peut lui
  supposer \'egalement toutes les valeurs depuis z\'ero jusqu'\`a
  l'unit\'e}'. See Jaynes \citey[\chap11]{jaynes1994_r1996} and Stigler
\citey{stigler1982} for interesting discussions and references on this
case. The case $\yk=1/2$, not discussed here, was advocated by
Jeffreys \citey{jeffreys1946}.

The third model $\ycom{\yl}$ in our comparison is defined by the suggestive
density
\begin{equation}
  \label{eq:G-mult}
  \yG(\yp \| \ycom{\yl})\,\ydp \defd \ync(\yl)
\frac{\yl!}{\prod_i (\yl p_i)!}\ydp,
\qquad \yl \ge 1;
\end{equation}
with $\ync(\yl)$ a normalization factor and, here and in the following,
\begin{equation}
  \label{eq:gamma_fact}
 x! \defd \Gamma(x + 1).
\end{equation}
I have been unable to characterize this model in more intuitive terms or to
derive it from particular assumptions as can be done for the Johnson one.
This is an interesting problem which deserves further study. Apart from the
normalization, the expression above is (intentionally, as we shall see
later) identical in form with the multiplicity factor~\eqref{eq:multfactor}
and for this reason I call this the \emph{multiplicity} model with
parameter $\yl$. We shall use the values $\yl=1, 2, 5, 50$ and $\yl$ large.
The case $\yl=1$ (a generalized density proportional to $\prod_i p_i^{-1}$)
was proposed by Haldane \citey{haldane1948} and is discussed by Jeffreys
\citey[\sect3.1, p.~123 ff.]{jeffreys1939_r2003} and Jaynes \citey[\sect
VII]{jaynes1968}; see also Zellner \citey[\sect2.13]{zellner1971_r1996}.


As you have already guessed, this model has been chosen because for
appropriate values of the parameter $\yl$ it gives the same distributions
as the \me\ principle.

Further mathematical properties of our models are discussed in
app.~\ref{cha:formulae}.

\chapter{Main results}
\label{cha:discussion}

The plausibilities assigned by the \me\ principle and our three
exchangeable models in the fifteen problems are derived
app.~\ref{cha:formulae} and presented in the tables of
app.~\ref{cha:tables}, p.~\pageref{cha:tables}. The features that strike me most in these tables are
the following.

\paragraph{Cases with one or two old throws, \ie\ $N=1\text{ or }2$, with an
  average $a=5$:}
In the first case we know for sure that the old throw must give \yC, so its
plausibility distribution must be $(0,0,0,0,100,0)\,\%$. In the second case
we know that the faces \yU, \yB, \yT\ cannot appear because the total sum
of the two old throws must be 10; so the distribution must have the form
$(0,0,0,\dotv,\dotv,\dotv)$. If we interpret the \me\ distribution as
referring to old throws, it is in these cases not just unreasonable, but
plainly \emph{wrong}. The \me\ principle is therefore not meant to be used
for old throws when $N$ is small.

All three exchangeable models give correct results instead. This was
expected: in situations of certainty the plausibility calculus reduces to
the truth-logical calculus.

\paragraph{Cases with $N=1\text{ or }2$ and $a=6$:} We know that the one or
two old throws give \yS. This is certainly no ground to suppose that the
throwing technique or the die be completely biased towards the face \yS; we
are in fact excluding any particular detailed background knowledge, as \eg\
that all faces of the die have the same, unknown, number of pips --- so
that in these cases that number is revealed to be 6. In general, knowledge
of the outcome of only one or two old throws should not make our
predictions for a new throw deviate very much from the uniform
distribution. If interpreted to refer to a new throw, the \me\ distribution
is therefore quite unreasonable since it concentrates all plausibility on
face \yS. The same is true for the $N=1\text{ or }2$, $a=5$ cases. The \me\
principle is therefore not meant to be applied for new throws when $N$ is
small.

The Johnson and multiplicity models on the other hand give
reasonably more uniform distributions for the new throw, especially for
larger values of the parameters $\yk$ and $\yl$. 

\bigskip

Note also how the exchangeable models respect the logical symmetries of
these cases: When $N=1$ the cases $a=6$ and $a=5$ are completely symmetric
under the exchange of those faces (of all faces, indeed). When $N=2$ and
$a=6$ our knowledge --- that only face \yS\ appears in old throws --- is
symmetric under exchange of all other faces; this symmetry is respected by
the exchangeable models in both old and new throws. When $N=2$ and $a=5$ we
know that the outcomes must be $(\yQ,\yS)$, $(\yC,\yC)$, or $(\yS,\yQ)$;
this symmetry under permutation of \yQ, \yS, and of \yU, \yB, \yT\ is
again respected in the exchangeable models.

\paragraph{Case with $N=1$ (or $N$ odd) and $a=7/2$:} An
obviously impossible case. This is reflected in the fact that plausibility
logic yields $\tfrac{0}{0}$, for old and new throws: any inference is
completely arbitrary and unreliable because we have been given
contradictory data. Yet the \me\ gives the uniform distribution.

\bigskip

We have concluded that the \me\ principle is not meant to be applied for
either old or new throws when $N$ is small. But how small is `small'?
Let us examine the examples with six old throws.

\paragraph{Case $N=6$, $a=6$:} All six old throws must have outcome \yS. The
\me\ principle and the exchangeable models give the correct distribution.
Knowing that six old throws give \yS, which plausibilities would you assign
for a new throw? I should be surprised at the set of sixes, but should not
conclude yet that the throwing technique or the die be completely biased
towards \yS, as the \me\ principle suggests instead. So I find the latter's
result unreasonable.

The multiplicity model with $\yl=50$ gives in my opinion the most reasonable
distribution for a new throw. Also the Johnson model with $\yk$ somewhere
between $5$ and $50$ gives a reasonable distribution.

\paragraph{Case $N=6$, $a=5$:} The \me\ distribution is reasonable when
interpreted as a distribution for old throws, but I still prefer the
Johnson and the multiplicity models with $\yk$ and $\yl$ around
$50$ that give a few percents less to the \yS\ and \yC\ faces.

For a new throw, the \me\ distribution is less reasonable; I find that
concentrating $75\,\%$ of the plausibility on \yC\ and \yS\ is too much.
Again, the multiplicity model with $\yl=50$ or the Johnson with
$\yk$ between $5$ and $50$ give results most reasonable for me.

\bigskip

So with an average based on six old throws the \me\ principle still gives
unsatisfactory answers. Let us look at the cases with twelve throws.

\paragraph{Case $N=12$, $a=6$:} We know again that all twelve old throws
give \yS; the \me\ principle and the three exchangeable models give the
correct answer. What about a new throw? I should start to believe that the
die or throwing technique be biased towards \yS, but should still not be
sure that they be \emph{completely} biased. So the \me\ principle's answer
is still unreasonable for me. I find the Johnson and the multiplicity
models with low parameter values, just above $1$, more reasonable, but
unreasonable with higher parameter values (\cf\ \sect\ref{cha:conclusions}).

\paragraph{Case $N=12$, $a=6\text{ or }5$:} The \me\ principle and all
exchangeable models give reasonable distributions for an old throw. For a
new throw I find the \me\ distribution too concentrated on the faces \yC\
and \yS; I prefer the exchangeable models with lower parameter values.

\bigskip

So for the problems with $N=12$ the \me\ principle gives more reasonable,
though not yet satisfactory, results.

Let us consider very large values of $N$ then. Here I find the \me\
distributions reasonable, for both old and new throws. How do the
exchangeable models compare?

\paragraph{Case with $N$ large, $a=6$:} The \me\ principle as well as the
Johnson and multiplicity models say that the larger the number of
old throws that give \yS, the larger the plausibility that the die or
throwing technique are completely biased toward face \yS\ for a new throw
as well; all plausibility is asymptotically concentrated on that face. This
is what I should indeed believe.

\paragraph{Case with $N$ large, $a=5$:} This case is the most interesting.
As usual, when the parameters $\yk$ or $\yl$ are much larger than $N$ both
the Johnson and multiplicity models behave like the fair-throw one, as
explained is app.~\ref{cha:formulae}. For parameter values which are large,
but small compared to the number of old throws, the Johnson model gives a
distribution asymptotically equal to that of the \me\ principle with the
\emph{Burg} entropy \citep{burg1975}, rather than the Shannon one. 

On the other hand, \emph{the multiplicity model gives exactly the same
  distribution as the usual maximum-\bd (Shannon-)\bd entropy principle}.
This result is extremely important and holds for both old and new throws.

\paragraph{Cases with $1<N\le 12$, $a=7/2$:} In all these
problems the \me\ principle and the exchangeable models for larger $\yk$
and $\yl$ give reasonable distributions, for both old and new throws. Note
how all exchangeable models including the fair-throw one give, for old
throws, slightly larger plausibilities to outcomes nearer \yT\ or \yQ.
This is not strange, as some counting shows. For example, among the 146
sets of four-throw outcomes summing up to 14 the face \yU\ appears 21 times
as first-throw outcome whereas the face \yT\ 27 times.

%

\bigskip

Ironically, the results of the exchangeable models have in many cases
greater Shannon entropies than those of the \me\ principle. If we regard
the Shannon entropy as a measure of `incertitude' in a plausibility
distribution, then the exchangeable models give in those cases more
`uncertain' or, as Jaynes would say, less `committing' answers than the
\me\ principle. This is not in contradiction with the principle, of course:
in those cases, the exchangeable models do not satisfy the constraint that
the average be equal to the expectation. See the discussion in the next
section.

\chapter{Conclusions and various remarks}
\label{cha:conclusions}

Most of the conclusions I state here are personal, in the sense that they
are not only a matter of logic but of taste as well. I invite you to
peruse the numerical results presented in the tables of
app.~\ref{cha:tables} and arrive at your own conclusions.

When the number of throws on which the average is calculated is small, the
\me\ principle gives very unreasonable or even wrong results, depending on
whether its distribution is interpreted as referring to new or old throws.
The performance gets better the larger the numbers of throws, but with
twelve of them I still find the \me\ distributions unreasonable.

On which grounds do I say `reasonable' or `unreasonable'? This is a very
important and interesting question. At first, my judgements were almost
instinctive; they came from background knowledge not immediately present to
the mind (see Jeffreys' colourful discussion on this point \citey[\sect3.1,
pp.~123--124]{jeffreys1939_r2003}). With some introspection, I can say my
grounds are these: I think it is very difficult to master such a throwing
technique or to construct such a regularly numbered die as would lead me to
assign a plausibility distribution remarkably different from the uniform
one. Therefore I choose, on the simplex of distributions, a density very
peaked around the uniform distribution. It is easy, on the other hand, to
construct a die with two or more faces showing the same number of pips. I
should reflect this by superposing to my previous density another
concentrated on the \emph{facets} of the simplex. Thus, neither the Johnson
nor the multiplicity model represents my background knowledge exactly. I
can also add that if the `die' were very irregular, \eg\ with one length
twice the other two, I should complain of having been given false data,
since I should not call that a `die'.

The plausibility model I choose reflects my knowledge about dice-throwing.
In a problem regarding something else, \eg\ the energy of some physical
system, my knowledge and hence the model used would probably be different. In
some situations it would even be reasonable to use non-exchangeable models.

The numerical comparison presented here shows that the use of plausibility
logic gives us more possibilities of `fine-tuning', of better representing
our background knowledges, than the principle of maximum entropy.
Enthusiasts for this principle would probably argue that in its full
generality it allows for finer tuning too: we can choose any convex region
on the simplex of distributions as constraint, representing the posterior
distributions we deem acceptable. But the full use of plausibility logic is
even more flexible: first, we can give different plausibilities to
different regions; second, the latter need not be convex. And, most
importantly, plausibility logic does not ask us to choose amongst
\emph{posterior} distributions, but to carefully specify a \emph{prior} one
on an appropriate hypothesis space --- it requires us to examine
\emph{whence} we start (including what question we are asking), not
\emph{whither} we want to arrive. And in inference problems this is always
advisable, lest we let our wishes, instead of the facts, suggest what is
more or less plausible.

Advocates of the \me\ principle may also contend that in the cases in
which this gives apparently wrong results (small $N$), it is because we have
chosen wrong or insufficient constraints. For instance, if we know that for
$N=2$ and $a=5$ it is logically impossible that some old throw give \yU,
\yB, or \yT, then we ought to impose this as a constraint. I might accept
this argument but still think that plausibility logic is superior since it
\emph{reveals} these constraints to us, as logical consequences of the
situation, without requiring us to put them again explicitly into the
theory. When I first obtained the results for the $N=2$, $a=5$ case I was
indeed surprised seeing that all exchangeable models give distributions of
the form $(0,0,0, x,y,x)$ for old throws and $(z,z,z, x,y,x)$ for a new
one. An analysis of the possible outcomes in this case then showed why this
must be logically so, as explained in \sect\ref{cha:discussion}. Simple
happenings like this show the beauty and power of plausibility logic.

As regards constraints, we have seen that in many cases the rule
`$\text{expectation} = \text{average}$' is not respected by the more
reasonable exchangeable models, which can therefore have \emph{greater}
entropy than the \me\ distribution. Because in some cases that constraint
rule is absurd. In general its range of applicability is `subjective', a
fact that is seldom stressed enough; consensus and thus objectivity are
only reached in limit cases. By `subjective' I do not mean `subject to
personal quirks', as de~Finetti seems to imply sometimes
\citey{definetti1931}, but `perceptibly dependent on almost imperceptible
differences in background knowledge' (as a chaotic dynamics on initial
conditions).

\bigskip

Although our study only concerns a special example, it is easy to see how
it can be generalized to a general theorem: The plausibility distribution
given by the \me\ principle \emph{for new observations} is the same as that
given by a particular class of infinitely exchangeable models in a
well-defined limit case. This was known, apparently, --- see Skilling
\citey{skilling1989b}, Rodr\'\i guez
\citey{rodriguez1991,rodriguez2002,rodriguez2003}, and references therein
--- but I have never seen it said explicitly. I prefer not to say that the
principle is `derived' from plausibility logic, for the former's
formulation is based on different axioms and primitives than the latter's,
and these we have not derived. But the \emph{procedure} to obtain the \me\
distribution can be justified by plausibility logic without the need of
those additional principles.

We need a characterization of the class of models from which the principle
stems, though. Interesting studies by Skilling, Rodr\'\i guez, Caticha
\amp\ Preuss try to characterize a model in this class by invoking the \me\
principle again; more about this below. It would be interesting to find a
characterization of the multiplicity model similar to that, mentioned in
\sect\ref{cha:models}, of the Johnson model; or in terms of a symmetry on a
particular hypothesis space. This characterization is also important
because one could try to generalize it to other hypothesis spaces beyond
the exchangeable-\bd model one (`plausibilities of plausibility-indexed
circumstances'). In this way we could obtain, from plausibility logic, the
form of an entropy for general statistical models (in the sense of Mielnik
and Holevo \citep{mielnik1974,mielnik1981,holevo1980_t1982,holevo1985});
see \eg\ the studies by Band \amp\ Park
\citey{bandetal1976,parketal1976,bandetal1977,parketal1977}, Slater
\citey{slater1991,slater1992,slater1992b,slater1992c,slater1992d}, Porta
Mana \amp\ Bj\"ork \citey{portamanaetal2005}, Barnum \etal\
\citey{barnumetal2009}.


The derivation of the \me\ procedure from plausibility logic is also useful
because it clearly shows in which situations the \me\ principle can be
applied, and pro\-vides reasonable results in those situations in which the
principle's answers are unreasonable or plainly wrong. The way this is
mathematically achieved is explained at the end of app.~\ref{cha:formulae}.
The derivation also makes clear that there may be situations in which we
can reasonably assign distributions different from those of the \me\
principle. For example, if we had reasons to use a Johnson model in our
inference, the conclusions of the \me\ principle (with Shannon's entropy)
would obviously be at variance with those
reasons.

The (large) parameter $\yl$ of the multiplicity model gives the order of
magnitude of the number of old throws $N$ at which the \me\ principle
begins to approximate our conclusions, as mathematically explained in
app.~\ref{cha:formulae}. How should the value of this parameter be chosen?
The answer depends on the problem, as does the choice of exchangeable ---
or non-exchangeable --- model. In our dice examples I should use a value of
$50$ or slightly larger. In a problem in which I could examine the die and
the throwing technique and get the impression that they are `fair', the
value would be even higher.

\bigskip

What place in inference has the \me\ principle then? As an `update rule' it
seems superfluous, since plausibility logic already provides us with such a
rule, which we have seen to give more satisfactory results. Is it a
procedure to assign `prior' plausibilities, as Jaynes continually stressed?
But then using average data would be inappropriate, because they are the
kind of data that can be used in Bayes' rule instead. Is it a principle for
selecting a prior distributions among those we deem appropriate? Perhaps:
we want a distribution that give highest plausibility to such-and-such
average; which to choose? Let's take that with maximum entropy. And yet, it
would be better to ask \emph{why} we want that such-and-such average have
highest plausibility. If is it because we have observed that average in
similar situations, why not just apply plausibility logic and Bayes' theorem,
as we have done in this study? Jaynes \citey[p.~27]{jaynes1988b} states
that the \me\ principle `is designed to cover more general situations,
where it does not make sense to speak of \textquotedbl
trials\textquotedbl'. But I have failed to find, even in Jaynes' writings,
any examples of such `more general situations' or at least situations not
involving some kind of repetitions of observations (`trials').

The mystery about the foundations of the \me\ principle remains. The reason
of the present study stemmed from my strictly personal opinion that any
`update rule' is \begin{inparaenum}[(a)]
\item\label{item:part_case} a special case of the general update rule of plausibility logic, or
\item\label{item:incons} inconsistent, or
\item\label{item:not} not an update rule;
\end{inparaenum}
and that any procedure for assigning prior plausibilities from some data is
\begin{inparaenum}[(a)]
\item a special case of a plausibility-logic updating within a particular
  model, or
\item inconsistent, or
\item not a procedure for assigning priors.
\end{inparaenum}
I wanted to be sure that the \me\ principle did not fall into the
(\ref{item:incons})~alternatives. My belief now is that the
(\ref{item:part_case})~alternatives hold; although there is a small
possibility that the~(\ref{item:not}) be right instead.

Any use of the \me\ principle involving observational data usurps the just
and enlightened throne of plausibility logic. Therefore the only place
suitable for the principle seems to be, not in the choice of a
distribution, but in the construction of prior \emph{densities} over a
space of such distributions, where no observational data are directly
involved. Here, however, we have an infinite-dimensional simplex and the
principle requires the prior definition of a `canonical' density
\citep[\sect4.b]{jaynes1963} (which is usually different from that
discussed in app.~\ref{cha:can_dens}): we have a chicken-and-egg problem.
And what kind of constraints should be chosen for the principle, in such a
space? The studies on `entropic priors' of Skilling
\citey{skilling1989b,skilling1998}, Rodr\'\i guez,
\citey{rodriguez1991,rodriguez1991b,rodriguez2002,rodriguez2003}, and
Caticha \amp\ Preuss \citep{caticha2001,catichaetal2004} are interesting in
this respect, although they still leave me largely unconvinced for the time
being, for reasons that may be explained elsewhere.

\bigskip

A final note: the above discussion concerns the logical place and exact
premises of the \me\ principle (which are necessary also for pedagogical
purposes), but does not affect its usefulness and efficacy, both witnessed
by the vast range of applications and the number of books written about,
and thanks to, this principle.


\begin{acknowledgements}
  Many thanks to Philip Goyal for encouragement and discussions; to R\"udiger
  Schack and Chris Fuchs for sharing their degrees of belief; to Dawn
  Bombay for her always prompt and smiling library help; to Thomas Hahn for
  his integration library \emph{Cuba}; to the developers and maintainers of
  \LaTeX, Emacs, AUC\TeX, MiK\TeX, arXiv; to Angela Hewitt for her
  marvellous interpretation of that Sarabande; and \emph{\langlatin{dulcis
      in fundo}} \langitalian{alle mie tre grazie per il loro supporto e
    continuo amore: Louise, Marianna, e soprattutto Miriam --- buon
    compleanno mami, sei la mamma pi\'u forte del mondo!}

My hate goes to PIP\reflectbox{P}O.

  Research at the Perimeter Institute is supported by the Government of Canada
  through Industry Canada and by the Province of Ontario through the
  Ministry of Research and Innovation.
\end{acknowledgements}

\newpage
\appendix

\newpage

\chapter{Tables of results}
\label{cha:tables}

Here are the distributions, for an old and a new throw, given by the \me\
principle and the fair-throw, Johnson, and multiplicity models described in
\sect\ref{cha:models}. The values are rounded to decimals of percentile;
this leads in many cases to unnormalized totals. The Shannon entropy of
each distribution is also given, within crotchets.

The general formul\ae\ used are derived in app.~\ref{cha:formulae}, where
mention is also made of the routines used, when necessary, for numerical
integration. See the same appendix also for explanation of the remarks
appearing under some distributions.










\newcommand{\dicetsize}{\footnotesize}
\newcommand{\numaver}[2]{\large #2, #1}

\newcommand{\uniformdist}{$\mathsf{(16.7,16.7,16.7,16.7,16.7,16.7)}$ {\scriptsize[$\mathsf{1.793}$]}}

\newcommand{\iuniformdist}{uniform distribution irrespective of $a$}

\newcommand{\commentd}[2]{\underset{\text{\vphantom{$\underline{a}$}{#2}}}{\underbrace{#1}}}

\newcommand{\commentdii}[1]{\commentd{#1}{like fair-throw model}}

\newcommand{\commentdme}[1]{\commentd{#1}{ME distribution}}

\newcommand{\commentdmb}[1]{\commentd{#1}{ME distribution for Burg entropy}}

\newcommand{\meuniformdist}{$\commentdme{\mathsf{(16.7,16.7,16.7,16.7,16.7,16.7)}}$ {\scriptsize[$\mathsf{1.793}$]}}

\newcommand{\mbuniformdist}{$\commentdmb{\mathsf{(16.7,16.7,16.7,16.7,16.7,16.7)}}$ {\scriptsize[$\mathsf{1.793}$]}}

\newcommand{\iiduniformdist}{$\commentdii{\mathsf{(16.7,16.7,16.7,16.7,16.7,16.7)}}$ {\scriptsize[$\mathsf{1.793}$]}}


\newcommand{\tablebt}
{
\centering{\dicetsize\sffamily
\begin{tabular}{lcc}
\multicolumn{3}{l}{\numaver{$a=\mathsf{7/2}$}{$N=\mathsf{2}$}}
\\
\toprule
model & old throw, $\p(R_i^1 \| \yA^N_a \land I)/\%$ {\scriptsize[$H/\mathsf{nat}$]}
& new throw, $\p(R^0_i \| \yA^N_a \land I)/\%$ {\scriptsize[$H/\mathsf{nat}$]} 
\\
\midrule
ME &  
\multicolumn{2}{c}{\uniformdist}
\\
\addlinespace
fair-t.\ $\yiid$ & $\mathsf{(16.7,16.7,16.7,16.7,16.7,16.7)}$ {\scriptsize[$\mathsf{1.793}$]}& 
\iuniformdist 
\\
\addlinespace
\multicolumn{2}{l}{Johnson $\ydir{\yk}$:} & \\
$\yk = \mathsf{1}$ &
$\mathsf{(16.7,16.7,16.7,16.7,16.7,16.7)}$ {\scriptsize[$\mathsf{1.793}$]}& 
$\mathsf{(16.7,16.7,16.7,16.7,16.7,16.7)}$ {\scriptsize[$\mathsf{1.793}$]} 
\\
$\yk = \mathsf{5}$ &
$\mathsf{(16.7,16.7,16.7,16.7,16.7,16.7)}$ {\scriptsize[$\mathsf{1.793}$]}& 
$\mathsf{(16.7,16.7,16.7,16.7,16.7,16.7)}$ {\scriptsize[$\mathsf{1.793}$]} 
\\
$\yk = \mathsf{50}$ &
$\mathsf{(16.7,16.7,16.7,16.7,16.7,16.7)}$ {\scriptsize[$\mathsf{1.793}$]}& 
$\mathsf{(16.7,16.7,16.7,16.7,16.7,16.7)}$ {\scriptsize[$\mathsf{1.793}$]} 
\\
$\yk$ large &
$\commentdii{\mathsf{(16.7,16.7,16.7,16.7,16.7,16.7)}}$ {\scriptsize[$\mathsf{1.793}$]}& 
\iuniformdist 
\\
\addlinespace
\multicolumn{2}{l}{multiplicity $\ycom{\yl}$:} & \\
$\yl = \mathsf{1}$ &
$\mathsf{(16.7,16.7,16.7,16.7,16.7,16.7)}$ {\scriptsize[$\mathsf{1.793}$]}& 
$\mathsf{(16.7,16.7,16.7,16.7,16.7,16.7)}$ {\scriptsize[$\mathsf{1.793}$]} 
\\
$\yl = \mathsf{5}$ &
$\mathsf{(16.7,16.7,16.7,16.7,16.7,16.7)}$ {\scriptsize[$\mathsf{1.793}$]}& 
$\mathsf{(16.7,16.7,16.7,16.7,16.7,16.7)}$ {\scriptsize[$\mathsf{1.793}$]} 
\\
$\yl = \mathsf{50}$ &
$\mathsf{(16.7,16.7,16.7,16.7,16.7,16.7)}$ {\scriptsize[$\mathsf{1.793}$]}& 
$\mathsf{(16.7,16.7,16.7,16.7,16.7,16.7)}$ {\scriptsize[$\mathsf{1.793}$]} 
\\
$\yl$ large &
$\commentdii{\mathsf{(16.7,16.7,16.7,16.7,16.7,16.7)}}$ {\scriptsize[$\mathsf{1.793}$]}& 
\iuniformdist 
\\
\bottomrule
\end{tabular}
}
}

\newcommand{\tablebs}
{
{\dicetsize\sffamily
\begin{tabular}{lcc}
\multicolumn{3}{l}{\numaver{$a=\mathsf{6}$}{$N=\mathsf{2}$}}
\\
\toprule
model & old throw, $\p(R_i^1 \| \yA^N_a \land I)/\%$ {\scriptsize[$H/\mathsf{nat}$]}
& new throw, $\p(R^0_i \| \yA^N_a \land I)/\%$ {\scriptsize[$H/\mathsf{nat}$]} 
\\
\midrule
ME &  
\multicolumn{2}{c}{$\mathsf{(0,0,0,0,0,100)}$ {\scriptsize[$\mathsf{0}$]}}
\\
\addlinespace
fair-t.\ $\yiid$ & $\mathsf{(0,0,0,0,0,100)}$ {\scriptsize[$\mathsf{0}$]}& 
\iuniformdist 
\\
\addlinespace
\multicolumn{2}{l}{Johnson $\ydir{\yk}$:} & \\
$\yk = \mathsf{1}$ &
$\mathsf{(0,0,0,0,0,100)}$ {\scriptsize[$\mathsf{0}$]}& 
$\mathsf{(12.5, 12.5, 12.5, 12.5, 12.5, 37.5)}$ {\scriptsize[$\mathsf{1.667}$]} 
\\
$\yk = \mathsf{5}$ &
$\mathsf{(0,0,0,0,0,100)}$ {\scriptsize[$\mathsf{0}$]}& 
$\mathsf{(15.6, 15.6, 15.6, 15.6, 15.6, 21.9)}$ {\scriptsize[$\mathsf{1.782}$]} 
\\
$\yk = \mathsf{50}$ &
$\mathsf{(0,0,0,0,0,100)}$ {\scriptsize[$\mathsf{0}$]}& 
$\mathsf{(16.6, 16.6, 16.6, 16.6, 16.6, 17.2)}$ {\scriptsize[$\mathsf{1.793}$]} 
\\
$\yk$ large &
$\commentdii{\mathsf{(0,0,0,0,0,100)}}$ {\scriptsize[$\mathsf{0}$]}& 
\iuniformdist 
\\
\addlinespace
\multicolumn{2}{l}{multiplicity $\ycom{\yl}$:} & \\
$\yl = \mathsf{1}$ &
$\mathsf{(0,0,0,0,0,100)}$ {\scriptsize[$\mathsf{0}$]}& 
$\mathsf{(12.6, 12.6, 12.6, 12.6, 12.6, 36.9)}$ {\scriptsize[$\mathsf{1.672}$]} 
\\
$\yl = \mathsf{5}$ &
$\mathsf{(0,0,0,0,0,100)}$ {\scriptsize[$\mathsf{0}$]}& 
$\mathsf{(13.5, 13.5, 13.5, 13.5, 13.5, 32.3)}$ {\scriptsize[$\mathsf{1.717}$]} 
\\
$\yl = \mathsf{50}$ &
$\mathsf{(0,0,0,0,0,100)}$ {\scriptsize[$\mathsf{0}$]}& 
$\mathsf{(16.0, 16.0, 16.0, 16.0, 16.0, 19.8)}$ {\scriptsize[$\mathsf{1.787}$]} 
\\
$\yl$ large &
$\commentdii{\mathsf{(0,0,0,0,0,100)}}$ {\scriptsize[$\mathsf{0}$]}& 
\iuniformdist 
\\
\bottomrule
\end{tabular}
}
}

\newcommand{\tablebc}
{
\centering{\dicetsize\sffamily
\begin{tabular}{lcc}
\multicolumn{3}{l}{\numaver{$a=\mathsf{5}$}{$N=\mathsf{2}$}}
\\
\toprule
model & old throw, $\p(R_i^1 \| \yA^N_a \land I)/\%$ {\scriptsize[$H/\mathsf{nat}$]}
& new throw, $\p(R^0_i \| \yA^N_a \land I)/\%$ {\scriptsize[$H/\mathsf{nat}$]} 
\\
\midrule
ME &  
\multicolumn{2}{c}{$\mathsf{(2.1,3.9,7.2,13.6,25.5,47.8)}$ {\scriptsize[$\mathsf{1.370}$]}}
\\
\addlinespace
fair-t.\ $\yiid$ & $\mathsf{(0, 0, 0, 33.3, 33.3, 33.3)}$ {\scriptsize[$\mathsf{1.10}$]}& 
\iuniformdist 
\\
\addlinespace
\multicolumn{2}{l}{Johnson $\ydir{\yk}$:} & \\
$\yk = \mathsf{1}$ &
$\mathsf{(0, 0, 0, 25.0, 50.0, 25.0)}$ {\scriptsize[$\mathsf{1.04}$]}& 
$\mathsf{(12.5, 12.5, 12.5, 18.8, 25.0, 18.8)}$ {\scriptsize[$\mathsf{1.755}$]} 
\\
$\yk = \mathsf{5}$ &
$\mathsf{(0, 0, 0, 31.2, 37.5, 31.2)}$ {\scriptsize[$\mathsf{1.09}$]}& 
$\mathsf{(15.6, 15.6, 15.6, 17.6, 18.0, 17.6)}$ {\scriptsize[$\mathsf{1.790}$]} 
\\
$\yk = \mathsf{50}$ &
$\mathsf{(0, 0, 0, 33.1, 33.8, 33.1)}$ {\scriptsize[$\mathsf{1.10}$]}& 
$\mathsf{(16.6, 16.6, 16.6, 16.8, 16.8, 16.8)}$ {\scriptsize[$\mathsf{1.793}$]} 
\\
$\yk$ large &
$\commentdii{\mathsf{(0, 0, 0, 33.3, 33.3, 33.3)}}$ {\scriptsize[$\mathsf{1.10}$]}& 
\iuniformdist 
\\
\addlinespace
\multicolumn{2}{l}{multiplicity $\ycom{\yl}$:} & \\
$\yl = \mathsf{1}$ &
$\mathsf{(0, 0, 0, 25.4, 49.1, 25.4)}$ {\scriptsize[$\mathsf{1.05}$]}& 
$\mathsf{(12.6, 12.6, 12.6, 18.8, 24.5, 18.8)}$ {\scriptsize[$\mathsf{1.756}$]} 
\\
$\yl = \mathsf{5}$ &
$\mathsf{(0, 0, 0, 26.9, 46.1, 26.9)}$ {\scriptsize[$\mathsf{1.06}$]}& 
$\mathsf{(13.4, 13.4, 13.4, 18.8, 22.1, 18.8)}$ {\scriptsize[$\mathsf{1.770}$]} 
\\
$\yl = \mathsf{50}$ &
$\mathsf{(0, 0, 0, 32.0, 35.9, 32.0)}$ {\scriptsize[$\mathsf{1.10}$]}& 
$\mathsf{(16.0, 16.0, 16.0, 17.3, 17.4, 17.3)}$ {\scriptsize[$\mathsf{1.791}$]} 
\\
$\yl$ large &
$\commentdii{\mathsf{(0, 0, 0, 33.3, 33.3, 33.3)}}$ {\scriptsize[$\mathsf{1.10}$]}& 
\iuniformdist 
\\
\bottomrule
\end{tabular}
}
}

\newcommand{\tableus}
{
{\dicetsize\sffamily
\begin{tabular}{lcc}
\multicolumn{3}{l}{\numaver{$a=\mathsf{6}$}{$N=\mathsf{1}$}}
\\
\toprule
model & old throw, $\p(R_i^1 \| \yA^N_a \land I)/\%$ {\scriptsize[$H/\mathsf{nat}$]}
& new throw, $\p(R^0_i \| \yA^N_a \land I)/\%$ {\scriptsize[$H/\mathsf{nat}$]} 
\\
\midrule
ME &  
\multicolumn{2}{c}{$\mathsf{(0,0,0,0,0,100)}$ {\scriptsize[$\mathsf{0}$]}}
\\
\addlinespace
fair-t.\ $\yiid$ & $\mathsf{(0,0,0,0,0,100)}$ {\scriptsize[$\mathsf{0}$]}& 
\iuniformdist 
\\
\addlinespace
\multicolumn{2}{l}{Johnson $\ydir{\yk}$:} & \\
$\yk = \mathsf{1}$ &
$\mathsf{(0,0,0,0,0,100)}$ {\scriptsize[$\mathsf{0}$]}& 
$\mathsf{(14.3,14.3,14.3,14.3,14.3,28.6)}$ {\scriptsize[$\mathsf{1.749}$]} 
\\
$\yk = \mathsf{5}$ &
$\mathsf{(0,0,0,0,0,100)}$ {\scriptsize[$\mathsf{0}$]}& 
$\mathsf{(16.1,16.1,16.1,16.1,16.1,19.4)}$ {\scriptsize[$\mathsf{1.788}$]} 
\\
$\yk = \mathsf{50}$ &
$\mathsf{(0,0,0,0,0,100)}$ {\scriptsize[$\mathsf{0}$]}& 
$\mathsf{(16.6,16.6,16.6,16.6,16.6,16.9)}$ {\scriptsize[$\mathsf{1.791}$]} 
\\
$\yk$ large &
$\commentdii{\mathsf{(0,0,0,0,0,100)}}$ {\scriptsize[$\mathsf{0}$]}& 
\iuniformdist 
\\
\addlinespace
\multicolumn{2}{l}{multiplicity $\ycom{\yl}$:} & \\
$\yl = \mathsf{1}$ &
$\mathsf{(0,0,0,0,0,100)}$ {\scriptsize[$\mathsf{0}$]}& 
$\mathsf{(14.4,14.4,14.4,14.4,14.4,28.2)}$ {\scriptsize[$\mathsf{1.752}$]} 
\\
$\yl = \mathsf{5}$ &
$\mathsf{(0,0,0,0,0,100)}$ {\scriptsize[$\mathsf{0}$]}& 
$\mathsf{(14.9,14.9,14.9,14.9,14.9,25.6)}$ {\scriptsize[$\mathsf{1.767}$]} 
\\
$\yl = \mathsf{50}$ &
$\mathsf{(0,0,0,0,0,100)}$ {\scriptsize[$\mathsf{0}$]}& 
$\mathsf{(16.3,16.3,16.3,16.3,16.3,18.3)}$ {\scriptsize[$\mathsf{1.789}$]} 
\\
$\yl$ large &
$\commentdii{\mathsf{(0,0,0,0,0,100)}}$ {\scriptsize[$\mathsf{0}$]}& 
\iuniformdist 
\\
\bottomrule
\end{tabular}
}
}

\newcommand{\tabless}
{
\centering{\dicetsize\sffamily
\begin{tabular}{lcc}
\multicolumn{3}{l}{\numaver{$a=\mathsf{6}$}{$N=\mathsf{6}$}}
\\
\toprule
model & old throw, $\p(R_i^1 \| \yA^N_a \land I)/\%$ {\scriptsize[$H/\mathsf{nat}$]}
& new throw, $\p(R^0_i \| \yA^N_a \land I)/\%$ {\scriptsize[$H/\mathsf{nat}$]} 
\\
\midrule
ME &  
\multicolumn{2}{c}{$\mathsf{(0,0,0,0,0,100)}$ {\scriptsize[$\mathsf{0}$]}}
\\
\addlinespace
fair-t.\ $\yiid$ & $\mathsf{(0,0,0,0,0,100)}$ {\scriptsize[$\mathsf{0}$]}& 
\iuniformdist 
\\
\addlinespace
\multicolumn{2}{l}{Johnson $\ydir{\yk}$:} & \\
$\yk = \mathsf{1}$ &
$\mathsf{(0,0,0,0,0,100)}$ {\scriptsize[$\mathsf{0}$]}& 
$\mathsf{(8.3,8.3,8.3,8.3,8.3,58.3)}$ {\scriptsize[$\mathsf{1.347}$]} 
\\
$\yk = \mathsf{5}$ &
$\mathsf{(0,0,0,0,0,100)}$ {\scriptsize[$\mathsf{0}$]}& 
$\mathsf{(13.9,13.9,13.9,13.9,13.9,30.6)}$ {\scriptsize[$\mathsf{1.734}$]} 
\\
$\yk = \mathsf{50}$ &
$\mathsf{(0,0,0,0,0,100)}$ {\scriptsize[$\mathsf{0}$]}& 
$\mathsf{(16.3,16.3,16.3,16.3,16.3,18.3)}$ {\scriptsize[$\mathsf{1.789}$]} 
\\
$\yk$ large &
$\commentdii{\mathsf{(0,0,0,0,0,100)}}$ {\scriptsize[$\mathsf{0}$]}& 
\iuniformdist 
\\
\addlinespace
\multicolumn{2}{l}{multiplicity $\ycom{\yl}$:} & \\
$\yl = \mathsf{1}$ &
$\mathsf{(0,0,0,0,0,100)}$ {\scriptsize[$\mathsf{0}$]}& 
$\mathsf{(8.5,8.5,8.5,8.5,8.5,57.5)}$ {\scriptsize[$\mathsf{1.366}$]} 
\\
$\yl = \mathsf{5}$ &
$\mathsf{(0,0,0,0,0,100)}$ {\scriptsize[$\mathsf{0}$]}& 
$\mathsf{(10.0,10.0,10.0,10.0,10.0,50.0)}$ {\scriptsize[$\mathsf{1.498}$]} 
\\
$\yl = \mathsf{50}$ &
$\mathsf{(0,0,0,0,0,100)}$ {\scriptsize[$\mathsf{0}$]}& 
$\mathsf{(15.0,15.0,15.0,15.0,15.0,24.8)}$ {\scriptsize[$\mathsf{1.769}$]} 
\\
$\yl$ large &
$\commentdii{\mathsf{(0,0,0,0,0,100)}}$ {\scriptsize[$\mathsf{0}$]}& 
\iuniformdist 
\\
\bottomrule
\end{tabular}
}
}

\newcommand{\tableds}
{
\centering{\dicetsize\sffamily
\begin{tabular}{lcc}
\multicolumn{3}{l}{\numaver{$a=\mathsf{6}$}{$N=\mathsf{12}$}}
\\
\toprule
model & old throw, $\p(R_i^1 \| \yA^N_a \land I)/\%$ {\scriptsize[$H/\mathsf{nat}$]}
& new throw, $\p(R^0_i \| \yA^N_a \land I)/\%$ {\scriptsize[$H/\mathsf{nat}$]} 
\\
\midrule
ME &  
\multicolumn{2}{c}{$\mathsf{(0,0,0,0,0,100)}$ {\scriptsize[$\mathsf{0}$]}}
\\
\addlinespace
fair-t.\ $\yiid$ & $\mathsf{(0,0,0,0,0,100)}$ {\scriptsize[$\mathsf{0}$]}& 
\iuniformdist 
\\
\addlinespace
\multicolumn{2}{l}{Johnson $\ydir{\yk}$:} & \\
$\yk = \mathsf{1}$ &
$\mathsf{(0,0,0,0,0,100)}$ {\scriptsize[$\mathsf{0}$]}& 
$\mathsf{(5.6,5.6,5.6,5.6,5.6,72.2)}$ {\scriptsize[$\mathsf{1.042}$]} 
\\
$\yk = \mathsf{5}$ &
$\mathsf{(0,0,0,0,0,100)}$ {\scriptsize[$\mathsf{0}$]}& 
$\mathsf{(11.9,11.9,11.9,11.9,11.9,40.5)}$ {\scriptsize[$\mathsf{1.633}$]} 
\\
$\yk = \mathsf{50}$ &
$\mathsf{(0,0,0,0,0,100)}$ {\scriptsize[$\mathsf{0}$]}& 
$\mathsf{(16.0,16.0,16.0,16.0,16.0,19.9)}$ {\scriptsize[$\mathsf{1.787}$]} 
\\
$\yk$ large &
$\commentdii{\mathsf{(0,0,0,0,0,100)}}$ {\scriptsize[$\mathsf{0}$]}& 
\iuniformdist 
\\
\addlinespace
\multicolumn{2}{l}{multiplicity $\ycom{\yl}$:} & \\
$\yl = \mathsf{1}$ &
$\mathsf{(0,0,0,0,0,100)}$ {\scriptsize[$\mathsf{0}$]}& 
$\mathsf{(5.7,5.7,5.7,5.7,5.7,71.6)}$ {\scriptsize[$\mathsf{1.056}$]} 
\\
$\yl = \mathsf{5}$ &
$\mathsf{(0,0,0,0,0,100)}$ {\scriptsize[$\mathsf{0}$]}& 
$\mathsf{(7.0,7.0,7.0,7.0,7.0,64.9)}$ {\scriptsize[$\mathsf{1.211}$]} 
\\
$\yl = \mathsf{50}$ &
$\mathsf{(0,0,0,0,0,100)}$ {\scriptsize[$\mathsf{0}$]}& 
$\mathsf{(13.9,13.9,13.9,13.9,13.9,30.7)}$ {\scriptsize[$\mathsf{1.734}$]} 
\\
$\yl$ large &
$\commentdii{\mathsf{(0,0,0,0,0,100)}}$ {\scriptsize[$\mathsf{0}$]}& 
\iuniformdist 
\\
\bottomrule
\end{tabular}
}
}

\newcommand{\tablels}
{
\centering{\dicetsize\sffamily
\begin{tabular}{lcc}
\multicolumn{3}{l}{\numaver{$a=\mathsf{6}$}{$N$ large}}
\\
\toprule
model & old throw, $\p(R_i^1 \| \yA^N_a \land I)/\%$ {\scriptsize[$H/\mathsf{nat}$]}
& new throw, $\p(R^0_i \| \yA^N_a \land I)/\%$ {\scriptsize[$H/\mathsf{nat}$]} 
\\
\midrule
ME &  
\multicolumn{2}{c}{$\mathsf{(0,0,0,0,0,100)}$ {\scriptsize[$\mathsf{0}$]}}
\\
\addlinespace
fair-t.\ $\yiid$ & $\commentdme{\mathsf{(0,0,0,0,0,100)}}$ {\scriptsize[$\mathsf{0}$]}& 
\iuniformdist 
\\
\addlinespace
\multicolumn{2}{l}{Johnson $\ydir{\yk}$:} & \\
$\yk = \mathsf{1}$ &
$\mathsf{(0,0,0,0,0,100)}$ {\scriptsize[$\mathsf{0}$]}& 
$\mathsf{(0,0,0,0,0,100)}$ {\scriptsize[$\mathsf{0}$]} 
\\
$\yk = \mathsf{5}$ &
$\mathsf{(0,0,0,0,0,100)}$ {\scriptsize[$\mathsf{0}$]}& 
$\mathsf{(0,0,0,0,0,100)}$ {\scriptsize[$\mathsf{0}$]} 
\\
$\yk = \mathsf{50}$ &
$\mathsf{(0,0,0,0,0,100)}$ {\scriptsize[$\mathsf{0}$]}& 
$\mathsf{(0,0,0,0,0,100)}$ {\scriptsize[$\mathsf{0}$]} 
\\
$\yk$ large:& \\
$N/\yk$ small &
$\commentdii{\mathsf{(0,0,0,0,0,100)}}$ {\scriptsize[$\mathsf{0}$]}& 
\iuniformdist 
\\
$\yk$ large:& \\
$N/\yk$ large &
$\commentdmb{\mathsf{(0,0,0,0,0,100)}}$ {\scriptsize[$\mathsf{0}$]}& 
$\commentdmb{\mathsf{(0,0,0,0,0,100)}}$ {\scriptsize[$\mathsf{0}$]} 
\\
\addlinespace
\multicolumn{2}{l}{multiplicity $\ycom{\yl}$:} & \\
$\yl = \mathsf{1}$ &
$\mathsf{(0,0,0,0,0,100)}$ {\scriptsize[$\mathsf{0}$]}& 
$\mathsf{(0,0,0,0,0,100)}$ {\scriptsize[$\mathsf{0}$]} 
\\
$\yl = \mathsf{5}$ &
$\mathsf{(0,0,0,0,0,100)}$ {\scriptsize[$\mathsf{0}$]}& 
$\mathsf{(0,0,0,0,0,100)}$ {\scriptsize[$\mathsf{0}$]} 
\\
$\yl = \mathsf{50}$ &
$\mathsf{(0,0,0,0,0,100)}$ {\scriptsize[$\mathsf{0}$]}& 
$\mathsf{(0,0,0,0,0,100)}$ {\scriptsize[$\mathsf{0}$]} 
\\
$\yl$ large:& \\
$N/\yl$ small &
$\commentdii{\mathsf{(0,0,0,0,0,100)}}$ {\scriptsize[$\mathsf{0}$]}& 
\iuniformdist 
\\
$\yl$ large:& \\
$N/\yl$ large &
$\commentdme{\mathsf{(0,0,0,0,0,100)}}$ {\scriptsize[$\mathsf{0}$]}& 
$\commentdme{\mathsf{(0,0,0,0,0,100)}}$ {\scriptsize[$\mathsf{0}$]} 
\\
\bottomrule
\end{tabular}
}
}

\newcommand{\tableuc}
{
\centering{\dicetsize\sffamily
\begin{tabular}{lcc}
\multicolumn{3}{l}{\numaver{$a=\mathsf{5}$}{$N=\mathsf{1}$}}
\\
\toprule
model & old throw, $\p(R_i^1 \| \yA^N_a \land I)/\%$ {\scriptsize[$H/\mathsf{nat}$]}
& new throw, $\p(R^0_i \| \yA^N_a \land I)/\%$ {\scriptsize[$H/\mathsf{nat}$]} 
\\
\midrule
ME &  
\multicolumn{2}{c}{$\mathsf{(2.1,3.9,7.2,13.6,25.5,47.8)}$ {\scriptsize[$\mathsf{1.370}$]}}
\\
\addlinespace
fair-t.\ $\yiid$ & $\mathsf{(0,0,0,0,100,0)}$ {\scriptsize[$\mathsf{0}$]}& 
\iuniformdist 
\\
\addlinespace
\multicolumn{2}{l}{Johnson $\ydir{\yk}$:} & \\
$\yk = \mathsf{1}$ &
$\mathsf{(0,0,0,0,100,0)}$ {\scriptsize[$\mathsf{0}$]}& 
$\mathsf{(14.3,14.3,14.3,14.3,28.6,14.3)}$ {\scriptsize[$\mathsf{1.749}$]} 
\\
$\yk = \mathsf{5}$ &
$\mathsf{(0,0,0,0,100,0)}$ {\scriptsize[$\mathsf{0}$]}& 
$\mathsf{(16.1,16.1,16.1,16.1,19.4,16.1)}$ {\scriptsize[$\mathsf{1.788}$]} 
\\
$\yk = \mathsf{50}$ &
$\mathsf{(0,0,0,0,100,0)}$ {\scriptsize[$\mathsf{0}$]}& 
$\mathsf{(16.6,16.6,16.6,16.6,16.9,16.6)}$ {\scriptsize[$\mathsf{1.791}$]} 
\\
$\yk$ large &
$\commentdii{\mathsf{(0,0,0,0,100,0)}}$ {\scriptsize[$\mathsf{0}$]}& 
\iuniformdist 
\\
\addlinespace
\multicolumn{2}{l}{multiplicity $\ycom{\yl}$:} & \\
$\yl = \mathsf{1}$ &
$\mathsf{(0,0,0,0,100,0)}$ {\scriptsize[$\mathsf{0}$]}& 
$\mathsf{(14.4,14.4,14.4,14.4,28.2,14.4)}$ {\scriptsize[$\mathsf{1.752}$]} 
\\
$\yl = \mathsf{5}$ &
$\mathsf{(0,0,0,0,100,0)}$ {\scriptsize[$\mathsf{0}$]}& 
$\mathsf{(14.9,14.9,14.9,14.9,25.6,14.9)}$ {\scriptsize[$\mathsf{1.767}$]} 
\\
$\yl = \mathsf{50}$ &
$\mathsf{(0,0,0,0,100,0)}$ {\scriptsize[$\mathsf{0}$]}& 
$\mathsf{(16.3,16.3,16.3,16.3,18.3,16.3)}$ {\scriptsize[$\mathsf{1.789}$]} 
\\
$\yl$ large &
$\commentdii{\mathsf{(0,0,0,0,100,0)}}$ {\scriptsize[$\mathsf{0}$]}& 
\iuniformdist 
\\
\bottomrule
\end{tabular}
}
}

\newcommand{\tableut}
{
\centering{\dicetsize\sffamily
\begin{tabular}{lcc}
\multicolumn{3}{l}{\numaver{$a=\mathsf{7/2}$}{$N=\mathsf{1}$}}
\\
\toprule
model & old throw, $\p(R_i^1 \| \yA^N_a \land I)/\%$ {\scriptsize[$H/\mathsf{nat}$]}
& new throw, $\p(R^0_i \| \yA^N_a \land I)/\%$ {\scriptsize[$H/\mathsf{nat}$]} 
\\
\midrule
ME &  
\multicolumn{2}{c}{\uniformdist}
\\
\addlinespace
all exch.\ models &
undefined &
undefined
\\
\bottomrule
\end{tabular}
}
}

\newcommand{\tablesc}
{
\centering{\dicetsize\sffamily
\begin{tabular}{lcc}
\multicolumn{3}{l}{\numaver{$a=\mathsf{5}$}{$N=\mathsf{6}$}}
\\
\toprule
model & old throw, $\p(R_i^1 \| \yA^N_a \land I)/\%$ {\scriptsize[$H/\mathsf{nat}$]}
& new throw, $\p(R^0_i \| \yA^N_a \land I)/\%$ {\scriptsize[$H/\mathsf{nat}$]} 
\\
\midrule
ME &  
\multicolumn{2}{c}{$\mathsf{(2.1,3.9,7.2,13.6,25.5,47.8)}$ {\scriptsize[$\mathsf{1.370}$]}}
\\
\addlinespace
fair-t.\ $\yiid$ & $\mathsf{(1.1,3.3,7.7,15.4,27.6,45.0)}$ {\scriptsize[$\mathsf{1.362}$]}& 
\iuniformdist 
\\
\addlinespace
\multicolumn{2}{l}{Johnson $\ydir{\yk}$:} & \\
$\yk = \mathsf{1}$ &
$\mathsf{(1.7,3.3,6.7,13.3,31.7,43.3)}$ {\scriptsize[$\mathsf{1.358}$]}& 
$\mathsf{(9.2,10.0,11.7,15.0,24.2,30.0)}$ {\scriptsize[$\mathsf{1.690}$]} 
\\
$\yk = \mathsf{5}$ &
$\mathsf{(1.4,3.5,7.3,14.7,27.5,45.5)}$ {\scriptsize[$\mathsf{1.363}$]}& 
$\mathsf{(14.1,14.5,15.1,16.3,18.5,21.5)}$ {\scriptsize[$\mathsf{1.780}$]} 
\\
$\yk = \mathsf{50}$ &
$\mathsf{(1.1,3.3,7.6,15.3,27.6,45.0)}$ {\scriptsize[$\mathsf{1.360}$]}& 
$\mathsf{(16.4,16.4,16.5,16.6,16.9,17.2)}$ {\scriptsize[$\mathsf{1.792}$]} 
\\
$\yk$ large &
$\commentdii{\mathsf{(1.1,3.3,7.7,15.4,27.6,45.0)}}$ {\scriptsize[$\mathsf{1.362}$]}& 
\iuniformdist 
\\
\addlinespace
\multicolumn{2}{l}{multiplicity $\ycom{\yl}$:} & \\
$\yl = \mathsf{1}$ &
$\mathsf{(1.7,3.4,6.7,13.4,31.2,43.6)}$ {\scriptsize[$\mathsf{1.360}$]}& 
$\mathsf{(9.2,10.1,11.8,15.1,23.8,29.9)}$ {\scriptsize[$\mathsf{1.691}$]} 
\\
$\yl = \mathsf{5}$ &
$\mathsf{(1.6,3.5,6.9,14.1,29.0,44.9)}$ {\scriptsize[$\mathsf{1.363}$]}& 
$\mathsf{(10.2,11.1,12.6,15.7,21.8,28.5)}$ {\scriptsize[$\mathsf{1.718}$]} 
\\
$\yl = \mathsf{50}$ &
$\mathsf{(1.3,3.4,7.4,14.9,27.5,45.4)}$ {\scriptsize[$\mathsf{1.361}$]}& 
$\mathsf{(15.0,15.2,15.7,16.5,17.9,19.7)}$ {\scriptsize[$\mathsf{1.787}$]} 
\\
$\yl$ large &
$\commentdii{\mathsf{(1.1,3.3,7.7,15.4,27.6,45.0)}}$ {\scriptsize[$\mathsf{1.362}$]}& 
\iuniformdist 
\\
\bottomrule
\end{tabular}
}
}

\newcommand{\tabledc}
{
\centering{\dicetsize\sffamily
\begin{tabular}{lcc}
\multicolumn{3}{l}{\numaver{$a=\mathsf{5}$}{$N=\mathsf{12}$}}
\\
\toprule
model & old throw, $\p(R_i^1 \| \yA^N_a \land I)/\%$ {\scriptsize[$H/\mathsf{nat}$]}
& new throw, $\p(R^0_i \| \yA^N_a \land I)/\%$ {\scriptsize[$H/\mathsf{nat}$]} 
\\
\midrule
ME &  
\multicolumn{2}{c}{$\mathsf{(2.1,3.9,7.2,13.6,25.5,47.8)}$ {\scriptsize[$\mathsf{1.370}$]}}
\\
\addlinespace
fair-t.\ $\yiid$ & $\mathsf{(1.6,3.6,7.4,14.4,26.6,46.4)}$ {\scriptsize[$\mathsf{1.366}$]}& 
\iuniformdist 
\\
\addlinespace
\multicolumn{2}{l}{Johnson $\ydir{\yk}$:} & \\
$\yk = \mathsf{1}$ &
$\mathsf{(2.7,4.3,6.6,11.7,26.6,48.2)}$ {\scriptsize[$\mathsf{1.367}$]}& 
$\mathsf{(7.3,8.4,9.9,13.4,23.3,37.7)}$ {\scriptsize[$\mathsf{1.605}$]} 
\\
$\yk = \mathsf{5}$ &
$\mathsf{(2.2,4.0,7.2,13.2,25.2,48.3)}$ {\scriptsize[$\mathsf{1.368}$]}& 
$\mathsf{(12.5,13.0,14.0,15.7,19.1,25.7)}$ {\scriptsize[$\mathsf{1.756}$]} 
\\
$\yk = \mathsf{50}$ &
$\mathsf{(1.7,3.6,7.4,14.3,26.3,46.7)}$ {\scriptsize[$\mathsf{1.367}$]}& 
$\mathsf{(16.1,16.2,16.3,16.6,17.0,17.8)}$ {\scriptsize[$\mathsf{1.791}$]} 
\\
$\yk$ large &
$\commentdii{\mathsf{(1.6,3.6,7.4,14.4,26.6,46.4)}}$ {\scriptsize[$\mathsf{1.366}$]}& 
\iuniformdist 
\\
\addlinespace
\multicolumn{2}{l}{multiplicity $\ycom{\yl}$:} & \\
$\yl = \mathsf{1}$ &
$\mathsf{(2.6,4.2,6.6,11.8,26.4,48.3)}$ {\scriptsize[$\mathsf{1.363}$]}& 
$\mathsf{(7.4,8.5,10.0,13.5,23.1,37.5)}$ {\scriptsize[$\mathsf{1.609}$]} 
\\
$\yl = \mathsf{5}$ &
$\mathsf{(2.5,4.1,6.8,12.4,25.6,48.5)}$ {\scriptsize[$\mathsf{1.365}$]}& 
$\mathsf{(8.1,9.2,10.9,14.4,22.1,35.2)}$ {\scriptsize[$\mathsf{1.645}$]} 
\\
$\yl = \mathsf{50}$ &
$\mathsf{(1.9,3.8,7.3,13.7,25.8,47.5)}$ {\scriptsize[$\mathsf{1.366}$]}& 
$\mathsf{(13.7,14.2,14.9,16.3,18.6,22.4)}$ {\scriptsize[$\mathsf{1.777}$]} 
\\
$\yl$ large &
$\commentdii{\mathsf{(1.6,3.6,7.4,14.4,26.6,46.4)}}$ {\scriptsize[$\mathsf{1.366}$]}& 
\iuniformdist 
\\
\bottomrule
\end{tabular}
}
}

\newcommand{\tablelc}
{
\centering{\dicetsize\sffamily
\begin{tabular}{lcc}
\multicolumn{3}{l}{\numaver{$a=\mathsf{5}$}{$N$ large}}
\\
\toprule
model & old throw, $\p(R_i^1 \| \yA^N_a \land I)/\%$ {\scriptsize[$H/\mathsf{nat}$]}
& new throw, $\p(R^0_i \| \yA^N_a \land I)/\%$ {\scriptsize[$H/\mathsf{nat}$]} 
\\
\midrule
ME &  
\multicolumn{2}{c}{$\mathsf{(2.1,3.9,7.2,13.6,25.5,47.8)}$ {\scriptsize[$\mathsf{1.370}$]}}
\\
\addlinespace
fair-t.\ $\yiid$ & $\commentdme{\mathsf{(2.1,3.9,7.2,13.6,25.5,47.8)}}$ {\scriptsize[$\mathsf{1.370}$]}& 
\iuniformdist 
\\
\addlinespace
\multicolumn{2}{l}{Johnson $\ydir{\yk}$:} & \\
$\yk = \mathsf{1}$ &
$\mathsf{(4.0,5.0,6.7,10.0,20.0,54.3)}$ {\scriptsize[$\mathsf{1.343}$]}& 
$\mathsf{(4.0,5.0,6.7,10.0,20.0,54.3)}$ {\scriptsize[$\mathsf{1.343}$]} 
\\
$\yk = \mathsf{5}$ &
$\mathsf{(4.3,5.3,6.9,9.8,17.2,56.5)}$ {\scriptsize[$\mathsf{1.328}$]}& 
$\mathsf{(4.3,5.3,6.9,9.8,17.2,56.5)}$ {\scriptsize[$\mathsf{1.328}$]} 
\\
$\yk = \mathsf{50}$ &
$\mathsf{(4.3,5.3,6.9,9.8,16.7,56.9)}$ {\scriptsize[$\mathsf{1.323}$]}& 
$\mathsf{(4.3,5.3,6.9,9.8,16.7,56.9)}$ {\scriptsize[$\mathsf{1.323}$]} 
\\
$\yk$ large:& \\
$N/\yk$ small &
$\commentdii{\mathsf{(2.1,3.9,7.2,13.6,25.5,47.8)}}$ {\scriptsize[$\mathsf{1.370}$]}& 
\iuniformdist 
\\
$\yk$ large:& \\
$N/\yk$ large &
$\commentdmb{\mathsf{(4.4,5.3,6.9,9.8,16.7,57.0)}}$ {\scriptsize[$\mathsf{1.325}$]}& 
$\commentdmb{\mathsf{(4.4,5.3,6.9,9.8,16.7,57.0)}}$ {\scriptsize[$\mathsf{1.325}$]} 
\\
\addlinespace
\multicolumn{2}{l}{multiplicity $\ycom{\yl}$:} & \\
$\yl = \mathsf{1}$ &
$\mathsf{(4.0,5.0,6.7,10.1,20.2,54.1)}$ {\scriptsize[$\mathsf{1.347}$]}& 
$\mathsf{(4.0,5.0,6.7,10.1,20.2,54.1)}$ {\scriptsize[$\mathsf{1.347}$]} 
\\
$\yl = \mathsf{5}$ &
$\mathsf{(3.6,4.7,6.6,10.7,21.8,52.5)}$ {\scriptsize[$\mathsf{1.352}$]}& 
$\mathsf{(3.6,4.7,6.6,10.7,21.8,52.5)}$ {\scriptsize[$\mathsf{1.352}$]} 
\\
$\yl = \mathsf{50}$ &
$\mathsf{(2.3,3.8,7.0,13.4,25.5,48.0)}$ {\scriptsize[$\mathsf{1.367}$]}& 
$\mathsf{(2.3,3.8,7.0,13.4,25.5,48.0)}$ {\scriptsize[$\mathsf{1.367}$]} 
\\
$\yl$ large:& \\
$N/\yl$ small &
$\commentdii{\mathsf{(2.1,3.9,7.2,13.6,25.5,47.8)}}$ {\scriptsize[$\mathsf{1.370}$]}& 
\iuniformdist 
\\
$\yl$ large:& \\
$N/\yl$ large &
$\commentdme{\mathsf{(2.1,3.9,7.2,13.6,25.5,47.8)}}$ {\scriptsize[$\mathsf{1.370}$]}& 
$\commentdme{\mathsf{(2.1,3.9,7.2,13.6,25.5,47.8)}}$ {\scriptsize[$\mathsf{1.370}$]} 
\\
\bottomrule
\end{tabular}
}
}

\newcommand{\tablest}
{
\centering{\dicetsize\sffamily
\begin{tabular}{lcc}
\multicolumn{3}{l}{\numaver{$a=\mathsf{7/2}$}{$N=\mathsf{6}$}}
\\
\toprule
model & old throw, $\p(R_i^1 \| \yA^N_a \land I)/\%$ {\scriptsize[$H/\mathsf{nat}$]}
& new throw, $\p(R^0_i \| \yA^N_a \land I)/\%$ {\scriptsize[$H/\mathsf{nat}$]} 
\\
\midrule
ME &  
\multicolumn{2}{c}{\uniformdist}
\\
\addlinespace
fair-t.\ $\yiid$ & $\mathsf{(15.0,17.0,18.0,18.0,17.0,15.0)}$ {\scriptsize[$\mathsf{1.789}$]}& 
\iuniformdist 
\\
\addlinespace
\multicolumn{2}{l}{Johnson $\ydir{\yk}$:} & \\
$\yk = \mathsf{1}$ &
$\mathsf{(13.0,16.1,20.8,20.8,16.1,13.0)}$ {\scriptsize[$\mathsf{1.772}$]}& 
$\mathsf{(14.8,16.4,18.8,18.8,16.4,14.8)}$ {\scriptsize[$\mathsf{1.787}$]} 
\\
$\yk = \mathsf{5}$ &
$\mathsf{(14.5,16.9,18.6,18.6,16.9,14.5)}$ {\scriptsize[$\mathsf{1.787}$]}& 
$\mathsf{(16.3,16.7,17.0,17.0,16.7,16.3)}$ {\scriptsize[$\mathsf{1.792}$]} 
\\
$\yk = \mathsf{50}$ &
$\mathsf{(15.0,17.0,18.1,18.1,17.0,15.0)}$ {\scriptsize[$\mathsf{1.790}$]}& 
$\mathsf{(16.6,16.7,16.7,16.7,16.7,16.6)}$ {\scriptsize[$\mathsf{1.792}$]} 
\\
$\yk$ large &
$\commentdii{\mathsf{(15.0,17.0,18.0,18.0,17.0,15.0)}}$ {\scriptsize[$\mathsf{1.789}$]}& 
\iuniformdist 
\\
\addlinespace
\multicolumn{2}{l}{multiplicity $\ycom{\yl}$:} & \\
$\yl = \mathsf{1}$ &
$\mathsf{(13.1,16.2,20.7,20.7,16.2,13.1)}$ {\scriptsize[$\mathsf{1.774}$]}& 
$\mathsf{(14.9,16.4,18.7,18.7,16.4,15.0)}$ {\scriptsize[$\mathsf{1.788}$]} 
\\
$\yl = \mathsf{5}$ &
$\mathsf{(13.5,16.6,20.0,20.0,16.5,13.5)}$ {\scriptsize[$\mathsf{1.780}$]}& 
$\mathsf{(15.4,16.7,18.0,18.0,16.6,15.4)}$ {\scriptsize[$\mathsf{1.791}$]} 
\\
$\yl = \mathsf{50}$ &
$\mathsf{(14.7,17.0,18.3,18.3,17.0,14.7)}$ {\scriptsize[$\mathsf{1.788}$]}& 
$\mathsf{(16.5,16.7,16.8,16.8,16.7,16.5)}$ {\scriptsize[$\mathsf{1.792}$]} 
\\
$\yl$ large &
$\commentdii{\mathsf{(15.0,17.0,18.0,18.0,17.0,15.0)}}$ {\scriptsize[$\mathsf{1.789}$]}& 
\iuniformdist 
\\
\bottomrule
\end{tabular}
}
}

\newcommand{\tabledt}
{
\centering{\dicetsize\sffamily
\begin{tabular}{lcc}
\multicolumn{3}{l}{\numaver{$a=\mathsf{7/2}$}{$N=\mathsf{12}$}}
\\
\toprule
model & old throw, $\p(R_i^1 \| \yA^N_a \land I)/\%$ {\scriptsize[$H/\mathsf{nat}$]}
& new throw, $\p(R^0_i \| \yA^N_a \land I)/\%$ {\scriptsize[$H/\mathsf{nat}$]} 
\\
\midrule
ME &  
\multicolumn{2}{c}{\uniformdist}
\\
\addlinespace
fair-t.\ $\yiid$ & $\mathsf{(15.9,16.8,17.3,17.3,16.8,15.9)}$ {\scriptsize[$\mathsf{1.791}$]}& 
\iuniformdist 
\\
\addlinespace
\multicolumn{2}{l}{Johnson $\ydir{\yk}$:} & \\
$\yk = \mathsf{1}$ &
$\mathsf{(13.5,16.5,20.0,20.0,16.5,13.5)}$ {\scriptsize[$\mathsf{1.779}$]}& 
$\mathsf{(14.5,16.6,18.9,18.9,16.6,14.5)}$ {\scriptsize[$\mathsf{1.786}$]} 
\\
$\yk = \mathsf{5}$ &
$\mathsf{(15.3,16.9,17.9,17.9,16.9,15.3)}$ {\scriptsize[$\mathsf{1.791}$]}& 
$\mathsf{(16.3,16.7,17.0,17.0,16.7,16.3)}$ {\scriptsize[$\mathsf{1.792}$]} 
\\
$\yk = \mathsf{50}$ &
$\mathsf{(15.8,16.8,17.4,17.4,16.8,15.8)}$ {\scriptsize[$\mathsf{1.791}$]}& 
$\mathsf{(16.6,16.7,16.7,16.7,16.7,16.6)}$ {\scriptsize[$\mathsf{1.792}$]} 
\\
$\yk$ large &
$\commentdii{\mathsf{(15.9,16.8,17.3,17.3,16.8,15.9)}}$ {\scriptsize[$\mathsf{1.791}$]}& 
\iuniformdist 
\\
\addlinespace
\multicolumn{2}{l}{multiplicity $\ycom{\yl}$:} & \\
$\yl = \mathsf{1}$ &
$\mathsf{(13.5,16.6,19.9,19.9,16.6,13.6)}$ {\scriptsize[$\mathsf{1.780}$]}& 
$\mathsf{(14.6,16.6,18.7,18.7,16.6,14.7)}$ {\scriptsize[$\mathsf{1.786}$]} 
\\
$\yl = \mathsf{5}$ &
$\mathsf{(14.1,16.8,19.1,19.1,16.8,14.1)}$ {\scriptsize[$\mathsf{1.784}$]}& 
$\mathsf{(15.2,16.8,18.0,18.0,16.7,15.2)}$ {\scriptsize[$\mathsf{1.789}$]} 
\\
$\yl = \mathsf{50}$ &
$\mathsf{(15.5,16.9,17.6,17.6,16.9,15.5)}$ {\scriptsize[$\mathsf{1.790}$]}& 
$\mathsf{(16.5,16.7,16.8,16.8,16.7,16.5)}$ {\scriptsize[$\mathsf{1.792}$]} 
\\
$\yl$ large &
$\commentdii{\mathsf{(15.9,16.8,17.3,17.3,16.8,15.9)}}$ {\scriptsize[$\mathsf{1.791}$]}& 
\iuniformdist 
\\
\bottomrule
\end{tabular}
}
}

\newcommand{\tablelt}
{
\centering{\dicetsize\sffamily
\begin{tabular}{lcc@{}}
\multicolumn{3}{l}{\numaver{$a=\mathsf{3.5}$}{$N$ large}}
\\
\toprule
model & old throw, $\p(R_i^1 \| \yA^N_a \land I)/\%$ {\scriptsize[$H/\mathsf{nat}$]}
& new throw, $\p(R^0_i \| \yA^N_a \land I)/\%$ {\scriptsize[$H/\mathsf{nat}$]} 
\\
\midrule
ME &  
\multicolumn{2}{c}{\uniformdist}
\\
\addlinespace
fair-t.\ $\yiid$ & \meuniformdist& 
\iuniformdist 
\\
\addlinespace
\multicolumn{2}{l}{Johnson $\ydir{\yk}$:} & \\
$\yk = \mathsf{1}$ &
$\mathsf{(14.1,16.6,19.3,19.3,16.6,14.1)}$ {\scriptsize[$\mathsf{1.784}$]}& 
$\mathsf{(14.1,16.6,19.3,19.3,16.6,14.1)}$ {\scriptsize[$\mathsf{1.784}$]} 
\\
$\yk = \mathsf{5}$ &
$\mathsf{(16.1,16.8,17.2,17.2,16.8,16.1)}$ {\scriptsize[$\mathsf{1.793}$]}& 
$\mathsf{(16.1,16.8,17.2,17.2,16.8,16.1)}$ {\scriptsize[$\mathsf{1.793}$]} 
\\
$\yk = \mathsf{50}$ &
$\mathsf{(16.6,16.7,16.7,16.7,16.7,16.6)}$ {\scriptsize[$\mathsf{1.792}$]}& 
$\mathsf{(16.6,16.7,16.7,16.7,16.7,16.6)}$ {\scriptsize[$\mathsf{1.792}$]} 
\\
$\yk$ large:& \\
$N/\yk$ small &
\iiduniformdist& 
\iuniformdist 
\\ 
$N/\yk$ large &
\mbuniformdist& 
\mbuniformdist 
\\
\addlinespace
\multicolumn{2}{l}{multiplicity $\ycom{\yl}$:} & \\
$\yl = \mathsf{1}$ &
$\mathsf{(14.2,16.6,19.2,19.2,16.6,14.2)}$ {\scriptsize[$\mathsf{1.784}$]}& 
$\mathsf{(14.2,16.6,19.2,19.2,16.6,14.2)}$ {\scriptsize[$\mathsf{1.784}$]} 
\\
$\yl = \mathsf{5}$ &
$\mathsf{(14.9,16.8,18.3,18.3,16.8,14.9)}$ {\scriptsize[$\mathsf{1.788}$]}& 
$\mathsf{(14.9,16.8,18.3,18.3,16.8,14.9)}$ {\scriptsize[$\mathsf{1.788}$]} 
\\
$\yl = \mathsf{50}$ &
$\mathsf{(16.5,16.7,16.8,16.8,16.7,16.5)}$ {\scriptsize[$\mathsf{1.792}$]}& 
$\mathsf{(16.5,16.7,16.8,16.8,16.7,16.5)}$ {\scriptsize[$\mathsf{1.792}$]} 
\\
$\yl$ large:& \\
$N/\yl$ small &
\iiduniformdist& 
\iuniformdist 
\\
$N/\yl$ large &
\meuniformdist& 
\meuniformdist 
\\
\bottomrule
\end{tabular}
}
}



{\centering\vspace*{\stretch{1}}
\tableus\\
\vspace*{\stretch{2}}
\tableuc\\
\vspace*{\stretch{1}}
\pagebreak[5]
\vspace*{\stretch{1}}
\tableut\\
\vspace*{\stretch{2}}


\tablebs\\
\vspace*{\stretch{2}}
\tablebc\\
\vspace*{\stretch{1}}
\pagebreak[5]

\vspace*{\stretch{1}}
\tablebt\\
\vspace*{\stretch{2}}
\tabless\\
\vspace*{\stretch{1}}

\pagebreak[5]

\vspace*{\stretch{1}}
\tablesc\\
\vspace*{\stretch{2}}
\tablest\\
\vspace*{\stretch{1}}

\pagebreak[5]
\vspace*{\stretch{1}}
\tableds\\
\vspace*{\stretch{2}}
\tabledc\\
\vspace*{\stretch{1}}

\pagebreak[5]

\vspace*{\stretch{1}}
\tabledt\\
\vspace*{\stretch{2}}
\tablels\\
\vspace*{\stretch{1}}

\pagebreak[5]

\vspace*{\stretch{1}}
\tablelc\\
\vspace*{\stretch{1}}

\pagebreak[5]
\vspace*{\stretch{1}}
\tablelt\\
\vspace*{\stretch{1}}
}

\newpage

\chapter{Formul{\ae}}
\label{cha:formulae}

In \sect\ref{cha:comparison} we introduced the propositions $\yR{j}{i}$
stating that throw $j$ shows face `$i$'. Throw $j=1$ is an `old' throw,
$j=0$ a `new' one, in the sense already explained in the same section. The
proposition $I$ denotes the model used in our inferences and other
background knowledge.

Let $F_{\yn}$ denote the statement that the number of occurrences of the
six possible outcomes in $N$ throws is $\yn \equiv (N_i)$. So $F_{\yn}$ is
a disjunction of conjunctions of $R$s; \eg, for $N=3$ and
$\yn=(0,0,0,2,1,0)$ (two \yQ\ and one \yC),
\begin{equation}
  \label{eq:F_eg}
F_{(0,0,0,2,1,0)} \equiv 
(\yR{1}{4} \land \yR{2}{4} \land \yR{3}{5}) \lor
(\yR{1}{4} \land \yR{2}{5} \land \yR{3}{4}) \lor
(\yR{1}{5} \land \yR{2}{4} \land \yR{3}{4}).
\end{equation}
With a finitely or infinitely exchangeable model $I$ for the old throws,
the plausibility of a set of outcomes depends on their frequencies but not
on which throws they occur; in our example,
\begin{equation}
  \label{eq:exch_freq_out_eg}
  \p(\yR{1}{4} \land \yR{2}{4} \land \yR{3}{5} \| I)
=
  \p(\yR{1}{4} \land \yR{2}{5} \land \yR{3}{4} \| I)
=
  \p(\yR{1}{5} \land \yR{2}{4} \land \yR{3}{5} \| I).
\end{equation}
In this case it is a simple combinatorial exercise to show that for any old
throw, \eg\ the first,
\begin{equation}
  \label{eq:RoldgivenF}
  \p(\yR{1}{i} \|F_{\yn} \& I) = N_i/N,
\end{equation}
unless $F_{\yn}$ and $I$ be incompatible, \ie\ $\p(F_{\yn}\|I)=0$, a case
which we exclude.

If the model is infinitely exchangeable with density $\yG(\yp\|I)\,\ydp$ we
have by standard combinatorial arguments
\citep(\eg)()[\sect1.2]{csiszaretal1981}[\sect4.3.2]{bernardoetal1994}[\sect2]{csiszaretal2004b}
\begin{equation}
  \label{eq:freq_from_exch}
  \p(F_{\yn}\|I) = \int 
N! \Biggl(\prod_l \frac{p_l^{N_l}}{N_l!}\Biggr)
\,\yG(\yp\|I)\,\ydp.
\end{equation}
Knowledge of the frequency of old throws leads to the `updated' density
\begin{equation}
  \label{eq:g_given_F}
  \yG(\yp\|F_{\yn} \land I)\,\ydp =
\frac{N! \biggl(\prod_l \frac{p_l^{N_l}}{N_l!}\biggr)
\,\yG(\yp\|I)}
{\p(F_{\yn}\|I)}
\ydp,
\end{equation}
from which we obtain the plausibility distribution for a new throw:
\begin{equation}
  \label{eq:RnewgivenF}
  \p(\yR{0}{i} \|F_{\yn} \& I) =
\int p_i \,\yG(\yp\|F_{\yn} \land I)\,\ydp
=
\frac{\int N! p_i \biggl(\prod_l \frac{p_l^{N_l}}{N_l!}\biggr)
\,\yG(\yp\|I)\,\ydp}
{\p(F_{\yn}\|I)}.
\end{equation}

Recall that $\yA^N_a$ is the statement that in $N$ throws the observed
average is $a$. This is obviously a statement about the possible outcome
frequencies of old throws, and we can write it as
\begin{equation}
\yo \inn \yn= aN, 
\qquad\text{with }
\yo \defd (1,2,3,4,5,6).
\label{eq:constraint}
\end{equation}
Hence
\begin{equation}
  \label{eq:aver_as_freq}
  \yA^N_a \equiv
\Lor_{\yo \inn \yn = aN} F_{\yn};
\end{equation}
\eg, if the observed average in two throws is $5/2$,
\begin{equation}
  \label{eq:eg_A}
  \yA^{(2)}_{5/2} \equiv
F_{(1,0,0,1,0,0)} \lor F_{(0,1,1,0,0,0)} \lor F_{(0,0,0,0,2,0)}.
\end{equation}
Sums over frequencies constrained by the formula $\yo \inn
\yn= aN$ will for brevity be denoted by 
$\sum_{\yn}^{(a)}$.

From \eqn\eqref{eq:aver_as_freq}, any plausibility conditional on
$\yA^N_a$ can  be resolved into a weighted sum of plausibilities
conditional on the $F_{\yn}$:
\begin{equation}
  \label{eq:decomp_freq}
  \p(\dotv \|\yA^N_a \land I) =
\frac{\sum_{\yn}^{(a)} \p(\dotv \|F_{\yn} \land I)\, \p(F_{\yn}\|I)}
{\sum_{\yn}^{(a)} \p(F_{\yn}\|I)}.
\end{equation}

Using formul\ae~\eqref{eq:RoldgivenF}--\eqref{eq:RnewgivenF} and
\eqref{eq:decomp_freq} we find
\begin{subequations}
  \begin{align}
    \p(\yR{1}{j} \|\yA^N_a \land I) &= \frac{\int \sum_{\yn}^{(a)}
      \tfrac{N_j}{N} \, N! \biggl(\prod_l\frac{p_l^{N_l}}{N_l!}\biggr)
      \,\yG(\yp\|I)\, \ydp} {\int\sum_{\yn}^{(a)} N!
      \biggl(\prod_l\frac{p_l^{N_l}}{N_l!}\biggr) \,\yG(\yp\|I)\, \ydp},
    \\
    \p(\yR{0}{j} \|\yA^N_a \land I) &= \frac{\int \sum_{\yn}^{(a)} p_j \,
      N! \biggl(\prod_l\frac{p_l^{N_l}}{N_l!}\biggr) \,\yG(\yp\|I)\, \ydp}
    {\int\sum_{\yn}^{(a)} N! \biggl(\prod_l\frac{p_l^{N_l}}{N_l!}\biggr)
      \,\yG(\yp\|I)\, \ydp}.
  \end{align}
\end{subequations}
With these formul\ae\ we can compute the plausibilities required in our
study.

\bigskip

Integration of these expressions is straightforward for the fair-throw
model. We obtain
\begin{subequations}
  \begin{align}
    \label{eq:fair_exact}
    \p(\yR{1}{i} \|\yA^N_a \land \yiid) &= \frac{\sum_{\yn}^{(a)}
      \tfrac{N_i}{N}\, 6^{-N}\, \prod_l \frac{N!}{N_l!}} {\sum_{\yn}^{(a)}
      6^{-N}\, \prod_l \frac{N!}{N_l!}},
    \\
    \p(\yR{0}{i} \|\yA^N_a \land \yiid) &= \frac{\sum_{\yn}^{(a)}
      \tfrac{1}{6}\, 6^{-N}\, \prod_l \frac{N!}{N_l!}} {\sum_{\yn}^{(a)}
      6^{-N}\, \prod_l \frac{N!}{N_l!}} = \frac{1}{6}.
  \end{align}
\end{subequations}
Note how the plausibility distribution for a new throw is independent from
any knowledge about old --- or any other --- throws. This model, like any
other \iid\ one, does not allow to `learn from experience'.

By \eqn\eqref{eq:multfactor} for enough large $N$ the first plausibility
above is approximated by
\begin{equation}
    \label{eq:fair_Nlarge}
    \begin{split}
      \p(\yR{1}{i} \|\yA^N_a \land \yiid) &\asy \frac{\sum_{\yn}^{(a)}
        \yf_i\, \exp[N H(\yff)] \, \ydf} {\sum_{\yn}^{(a)} \exp[N H(\yff)]
        \, \ydf} \quad\text{with $\yf_i \defd N_i/N$}
      \\
      &\asy \yf_i \text{ which maximizes $H(\yff)$ under the constraint
        $\yo\inn\yff=a$},
    \end{split}
\end{equation}
where $H$ is Shannon's entropy \eqref{eq:shan-entropy}. This is the result
of van Campenhout and Cover \citey{vancampenhoutetal1981} and Csisz\'ar
\citey{csiszar1985}: for old throws, the \me\ distribution is the
as\-ymp\-tot\-ic distribution for an old throw conditional on the average of a
large number of old throws, under the assumption of a uniform \iid\ model.

\bigskip

In the case of the Johnson model, with density
\begin{equation}
  \labelbis{eq:G-dir}
  \yG(\yp \| \ydir{\yk})\,\ydp = 
\G(1+\lsum_l\yk)\,\Biggl[\prod_l \frac{p_l^{\yk-1}}{\G(\yk)}\Biggr]\ydp,
\qquad \yk > 0,
\end{equation}
the integrals above can be reduced to Dirichlet's generalization of the beta
integral \citep[\sect1.8]{andrewsetal1999_r2000},
\begin{equation}
  \label{eq:beta_int}
  \int \bigl(\lprod_l p_l^{b_l-1}\bigr) \,\ydp =
  \frac{\prod_l \G(b_l)}{\G(1+\sum_lb_l)},
  \qquad \Re b_l > 0,
\end{equation}
and after some simplifications we obtain the closed-form
formul\ae
\begin{subequations}
  \begin{align}
    \label{eq:dirich_exact}
    \p(\yR{1}{i} \|\yA^N_a \land \ydir{\yk}) &= \frac{\sum_{\yn}^{(a)}
      \tfrac{N_i}{N} \prod_l \frac{(N_l+\yk-1)!}{N_l!}} {\sum_{\yn}^{(a)}
      \prod_i \frac{(N_l+\yk-1)!}{N_l!}},
    \\
    \p(\yR{0}{i} \|\yA^N_a \land \ydir{\yk}) &= \frac{\sum_{\yn}^{(a)}
      \tfrac{N_i+\yk}{N+6\yk} \prod_l \frac{(N_l+\yk-1)!}{N_l!}}
    {\sum_{\yn}^{(a)} \prod_l \frac{(N_l+\yk-1)!}{N_l!}}.
  \end{align}
\end{subequations}

Note how old throws are weighted means of $N_i/N$, new throws of
$(N_i+\yk)/(N+6\yk)$. The Johnson model behaves as if we knew about the
existence of $6\yk$ additional old throws in which each face occurred $\yk$
times. \Cf\ Jaynes \citey{jaynes1986d_r1996}, Johnson \citey{johnson1932c},
Zabell \citey{zabell1982}.

If $\yk$ is large we can use \eqn\eqref{eq:multfactor} to show that the
above formul\ae\ asymptotically become
\begin{subequations}
  \begin{align}
    \label{eq:dirich_klarge}
    \p(\yR{1}{i} \|\yA^N_a \land \ydir{\yk}) &\asy \frac{\sum_{\yn}^{(a)}
      \tfrac{N_i}{N}\, 6^{-N}\, \prod_l \frac{N!}{N_l!}} {\sum_{\yn}^{(a)}
      6^{-N}\, \prod_l \frac{N!}{N_l!}},
    \\
    \p(\yR{0}{i} \|\yA^N_a \land \ydir{\yk}) &\asy \frac{\sum_{\yn}^{(a)}
      \tfrac{1}{6}\, 6^{-N}\, \prod_l \frac{N!}{N_l!}} {\sum_{\yn}^{(a)}
      6^{-N}\, \prod_l \frac{N!}{N_l!}} = \frac{1}{6},
  \end{align}
\end{subequations}
and the Johnson model is approximated by the fair-throw one. This is also
true for $N$ large but enough smaller than $\yk$.

When $N$ is enough large and $\yk$ finite \eqn\eqref{eq:multfactor} can
again be used to show that  expressions~\eqref{eq:dirich_exact} are both
approximated by the integral
\begin{equation}
  \label{eq:dirich_Nlarge}
  \p(\yR{1}{i} \|\yA^N_a \land \ydir{\yk})
\asy
  \p(\yR{0}{i} \|\yA^N_a \land \ydir{\yk})
\asy
\frac{\int_\ysimpa
\yf_i \,
\prod_l \yf_l^{\yk-1}
\,\ydf}
{\int_\ysimpa
\prod_l \yf_l^{\yk-1}
\,\ydf},
\end{equation}
where $\ysimpa$ is the intersection between the plausibility simplex
$\ysimp$~\eqref{eq:simplex} and the constraint hyperplane $\yo\inn\yff =a$.
If $\yk$ is also enough large but enough smaller than $N$ the integral
above becomes
\begin{equation}
  \label{eq:to_burg}
\frac{\int_\ysimpa
\yf_i \,
\prod_l \yf_l^{\yk-1}
\,\ydf}
{\int_\ysimpa
\prod_l \yf_l^{\yk-1}
\,\ydf}
\asy
\frac{\int_\ysimpa
\yf_i \,
\exp[\yk \ybu(\yff)]
\,\ydf}
{\int_\ysimpa
\exp[\yk \ybu(\yff)]
\,\ydf}
\end{equation}
where $\ybu$ is Burg's entropy \citep{burg1975}:
\begin{equation}
  \label{eq:Burg_ent}
  \ybu(\yff) \defd 
\lsum_i \ln \yf_i.
\end{equation}
Then
\begin{multline}
  \label{eq:dirich_Nlarge_Klarge}
  \p(\yR{1}{i} \|\yA^N_a \land \ydir{\yk})
\asy
  \p(\yR{0}{i} \|\yA^N_a \land \ydir{\yk})
\asy{}
\\
\yf_i
\text{ which maximizes $\ybu(\yff)$ under the constraint $\yo\inn\yff=a$}.
\end{multline}
Thus the Johnson model yields the \me\ principle with \emph{Burg}'s
entropy; \cf\ Jaynes \citey[p.~19]{jaynes1986d_r1996}. Had we used a
generalized Johnson model with density
\begin{equation}
  \label{eq:G-dir_gen}
  \yG(\yp \| \ydir{\yk,\ykkk})\,\ydp = 
\G(1+\yk)\,\Biggl[\prod_l \frac{p_l^{\yk\ykk_l-1}}{\G(\yk\ykk_l)}\Biggr]\ydp,
\qquad \text{$\yk,\ykk_i> 0$, $\lsum_i\ykk_i =1$}
\end{equation}
the asymptotic expression~\eqref{eq:dirich_Nlarge_Klarge} would have had
the Kullback-Leibler divergence
\citep{kullbacketal1951,kullback1959_r1978,hobsonetal1973}
\begin{equation}
-\yD(\ykkk,\yff) \defd -\tsum_i \ykk_i \ln(\ykk_i/\yf_i)
\label{eq:KL}
\end{equation}
in place of Burg's
entropy $\ybu(\yff)$.

For the case in which both $N$ and $\yk$ are large and comparable with each
other, see the following analogous discussion for the multiplicity model.

\bigskip

The integrals for the multiplicity model, with density
\begin{equation}
  \labelbis{eq:G-mult}
  \yG(\yp \| \ycom{\yl})\,\ydp = \ync(\yl)
\frac{\yl!}{\prod_i (\yl p_i)!}\ydp,
%
\qquad \yl \ge 1,
\end{equation}
are
\begin{subequations}\label{eq:multipl}
  \begin{align}
    \p(\yR{1}{i} \|\yA^N_a \land \ycom{\yl}) &= \frac{\int \sum_{\yn}^{(a)}
      \tfrac{N_i}{N} \, N! L! \biggl(\prod_l\frac{p_l^{N_l}}{N_l! (\yl
        p_l)!}\biggr) \,\ydp} {\int\sum_{\yn}^{(a)} N! L!
      \biggl(\prod_l\frac{p_l^{N_l}}{N_l! (\yl p_l)!}\biggr) \,\ydp},
    \\
    \p(\yR{0}{i} \|\yA^N_a \land \ycom{\yl}) &= \frac{\int \sum_{\yn}^{(a)}
      p_i \, N! L! \biggl(\prod_l\frac{p_l^{N_l}}{N_l! (\yl p_l)!}\biggr)
      \,\ydp} {\int\sum_{\yn}^{(a)} N! L!
      \biggl(\prod_l\frac{p_l^{N_l}}{N_l! (\yl p_l)!}\biggr) \,\ydp}.
  \end{align}
\end{subequations}

When $\yl$ is large enough, the usual relationship~\eqref{eq:multfactor}
between Shannon's entropy and the multiplicity factor leads to the
approximations
\begin{subequations}\label{eq:multipl_Llarge}
  \begin{align}
    \p(\yR{1}{i} \|\yA^N_a \land \ycom{\yl}) &\asy \frac{\int
      \sum_{\yn}^{(a)} \tfrac{N_i}{N} \, N!
      \biggl(\prod_l\frac{p_l^{N_l}}{N_l!}\biggr) \, \exp[\yl H(\yp)]
      \,\ydp} {\int\sum_{\yn}^{(a)} N!
      \biggl(\prod_l\frac{p_l^{N_l}}{N_l!}\biggr) \, \exp[\yl H(\yp)]
      \,\ydp}
    \asy \frac{\sum_{\yn}^{(a)} \tfrac{N_j}{N}\, 6^{-N}\, \prod_l
      \frac{N!}{N_l!}} {\sum_{\yn}^{(a)} 6^{-N}\, \prod_l \frac{N!}{N_l!}},
    \\
    \p(\yR{0}{i} \|\yA^N_a \land \ycom{\yl}) &\asy \frac{\int
      \sum_{\yn}^{(a)} p_i \, \, N!
      \biggl(\prod_l\frac{p_l^{N_l}}{N_l!}\biggr) \, \exp[\yl H(\yp)]
      \,\ydp} {\int\sum_{\yn}^{(a)} N!
      \biggl(\prod_l\frac{p_l^{N_l}}{N_l!}\biggr) \, \exp[\yl H(\yp)]
      \,\ydp}
    \asy \frac{\sum_{\yn}^{(a)} \tfrac{1}{6}\, 6^{-N}\, \prod_l
      \frac{N!}{N_l!}} {\sum_{\yn}^{(a)} 6^{-N}\, \prod_l \frac{N!}{N_l!}}
    = \frac{1}{6};
  \end{align}
\end{subequations}
thus the multiplicity model is approximated by the fair-throw one for
enough large values of $\yl$, also when $N$ is large but smaller than $\yl$.

What happens if $N$ is large and $\yl$ finite? As for the Johnson model,
the integrals can be approximated by their restriction to the hyperplane
$\yo\inn\yff =a$:
\begin{equation}
  \label{eq:multip_Nlarge}
  \p(\yR{1}{i} \|\yA^N_a \land \ycom{\yl})
\asy
  \p(\yR{0}{i} \|\yA^N_a \land \ycom{\yl})
\asy
\frac{\int_\ysimpa
\yf_i \,
\bigl[\prod_l (\yl \yf_l)!\bigr]^{-1}
\ydf}
{\int_\ysimpa
\bigl[\prod_l (\yl \yf_l)!\bigr]^{-1}
\ydf}.
\end{equation}

Most interesting is the case in which $\yl$ is also large but still smaller
than $N$. From \eqn\eqref{eq:multfactor} as usual we see that the
integrals above are approximated by
\begin{equation}
  \label{eq:to_shannon}
\frac{\int_\ysimpa
\yf_i \,
\bigl[\prod_l (\yl \yf_l)!\bigr]^{-1}
\ydf}
{\int_\ysimpa
\bigl[\prod_l (\yl \yf_l)!\bigr]^{-1}
\ydf}
\asy
\frac{\int_\ysimpa
\yf_i \,
\exp[\yl H(\yff)]
\,\ydf}
{\int_\ysimpa
\exp[\yl H(\yff)]
\,\ydf}
\end{equation}
so that
\begin{multline}
  \label{eq:multip_Nlarge_Llarge}
  \p(\yR{1}{i} \|\yA^N_a \land \ycom{\yl})
\asy
  \p(\yR{0}{i} \|\yA^N_a \land \ycom{\yl})
\asy{}\\
\yf_i
\text{ which maximizes $H(\yff)$ under the constraint $\yo\inn\yff=a$}.
\end{multline}
That is, for large $N$, $\yl$, and $N/\yl$, the distribution
given by the multiplicity model for old and new throws is equal to that of
the \me\ principle.

It is easy to see why. For $N$ large, and larger than $\yl$, the data make
the posterior density of the model very peaked on the hyperplane determined
by the average $a$. On this hyperplane, however, the posterior density is
proportional to the prior one. For large values of the parameter $\yl$ the
density of the multiplicity model tends to have isopycnals (contour lines
of same density, with respect to the canonical density $\ydp$) coinciding
with the isentropes of the simplex of plausibility distributions, and very
peaked on those of larger entropy. Hence the final distribution corresponds
to that of larger entropy in the `constraint' hyperplane determined by the
average.

The assigned distribution is thus, in general, determined by the
competition be\-tween two peaks: that of the prior density, centred on the
uniform distribution, and that of the likelihood of the data, concentrated
on the `constraint' hyperplane. For a fixed, large $\yl$ and small $N$ the
first peak dominates and the assigned distribution is near the uniform one.
As $N$ increases the assigned distribution moves towards the constraint
hyperplane; and when $N$ becomes much larger than $\yl$ it is practically
on that hyperplane, though its exact position therein is still determined
by the prior density. This is why the multiplicity model gives reasonable
distributions for small $N$ but can approximate the \me\ distribution for
large $N$.

Of course this property is not exclusive to this model. Any other model
with isopycnals coinciding with isentropes and very peaked on those of
higher entropy will lead to the same distribution in problems with $N$
large enough; \eg\ a model with the `entropy prior' of Skilling, Rodr\'\i
guez, Caticha \amp\ Preuss discussed in \sect\ref{cha:conclusions}.

Using the generalized multiplicity model defined by the density
\begin{equation}
  \label{eq:G-mult-gen}
  \yG(\yp \| \ycom{\yl,\ykkk})\,\ydp = \ync(\yl,\ykkk)
\,\yl!\prod_i\frac{\ykk_i^{\yl p_i}}{(\yl p_i)!}\ydp,
%
\qquad\text{$\yl \ge 1$, $\ykk_i \ge 0$, $\lsum_i \ykk_i =1$}, 
\end{equation}
the asymptotic result~\eqref{eq:multip_Nlarge_Llarge} generalizes to the
`maximum-\bd relative-\bd entropy' principle,
\begin{multline}
  \label{eq:multip_Nlarge_Llarge}
  \p(\yR{1}{i} \|\yA^N_a \land \ycom{\yl,\ykkk})
\asy
  \p(\yR{0}{i} \|\yA^N_a \land \ycom{\yl,\ykkk})
\asy{}\\
\yf_i
\text{ which minimizes $\yD(\yff, \ykkk)$ 
under the constraint $\yo\inn\yff=a$},
\end{multline}
with the Kullback-Leibler divergence~\eqref{eq:KL} instead of Shannon's
entropy. Note that the roles of $\yff$
and $\ykkk$ are interchanged with respect to the Johnson model's case.

\bigskip

The integrals~\eqref{eq:dirich_Nlarge}, \eqref{eq:multipl},
\eqref{eq:multip_Nlarge} were numerically calculated with both the Monte
Carlo routine \emph{Suave} and the deterministic routine \emph{Cuhre} of
Hahn's multi\-di\-men\-sion\-al-integration library \emph{Cuba}
\citep{hahn2005}, comparing the results of both routines to appraise their
mutual consistency and precision. In most cases Suave was fastest and most
precise.

\chapter{The canonical density on a plausibility simplex: a whimsical definition}
\label{cha:can_dens}

Integrations over a plausibility simplex $\ysimp$ of dimension $n$ are
usually written as
\begin{equation}
  \label{eq:inte_nor}
  n! 
\int_0^1 \int_0^{1-p_1} \int_0^{1-p_1-p_2} \dotsi \int_0^{1- \sum_j^{n-2}p_j} 
f(p_i) \,\di p_{n-1} \dotsm \di p_3 \,\di p_2\,\di p_1
\end{equation}
or equivalently and more symmetrically as
\begin{equation}
  \label{eq:inte_delt}
  n! \int_0^1 \int_0^1 \dotsi \int_0^1 
g(p_i) \,\delt(1-\lsum p_i)
\,\di p_n \dotsm \di p_2\,\di p_1
\end{equation}
or even
\begin{equation}
  \label{eq:inte_delt}
  n! \int_0^\infty \int_0^\infty \dotsi \int_0^\infty 
g(p_i) \,\delt(1-\lsum p_i)
\,\di p_n \dotsm \di p_2\,\di p_1.
\end{equation}
When $g \equiv 1$ these integrals give the simplex a unit volume.

The first expression is unpleasant because asymmetric in the symmetric
variables $p_i$; the other two are unpleasant because they unnecessarily
invoke a generalized function. Behind these expressions there is a simple
and well-behaved canonical density or volume element $\yden$ over the
simplex, which gives the latter a unit volume. By density I mean a twisted,
positively oriented $n$-form. Twisted, or `odd', because its integration
does not require an inner orientation of the simplex (see Schouten
\citey{schoutenetal1940}[\chaps II, III]{schouten1951_r1989}, Burke
\citey{burke1983}[\chap IV]{burke1985_r1987}{burke1995}; also Marsden
\etal\ \citey[\chap7]{marsdenetal1983_r2002} and Choquet-Bruhat \etal\
\citey[\sect IV.B.1]{choquet-bruhatetal1977_r1996}); in fact, choosing an
inner orientation would break the permutation invariance of the functions
$p_i$.

This canonical density can be defined in two ways.

First way. Any $n$-dimensional convex set that can be affinely mapped onto
a finite region $\RR^n$ by a map $F$ can be given a canonical density
$\ydens$ by pulling back the the canonical density of
$\RR^n$ and rescaling:
\begin{equation}
  \label{eq:can_dens}
  \ydens \defd \frac{F^* \abs{\di\bm{x}}}{\int F^* \abs{\di\bm{x}}},
\end{equation}
where $\abs{\di\bm{x}}$ is the canonical density (twisted, positively
oriented $n$-form) of $\RR^n$. The rescaling gives the convex set a unit
volume. It is easily proven that this definition is independent of the
affine map $F$ chosen. When the convex set is a simplex, $\ydens$ is the
canonical density $\yden$.

The second way does not involve any embeddings, but uses instead the
naturally defined plausibility functions $p_i\colon \ysimp \to \RR$, which
characterize the simplex as a \emph{plausibility} simplex. The canonical
density is then implicitly defined as
\begin{equation}
  \label{eq:can_dens_impl}
  \yden \defd \sum_{\text{\makebox[3em]{$\set{i_1, \dotsc, i_{n-1}} 
\subset \set{1, \dotsc, n}$}}}
\abs{\di p_{i_1} \we \dotsb \we\di p_{i_{n-1}}},
\end{equation}
the indices running over all permutations of $n-1$ elements of $\set{1,
  \dotsc, n}$. The magnitude operator $\abs{\dotv}$ transforms any density
into the twisted, positively oriented equipollent one. The expression above
is symmetric on the $p_i$ and does not involve generalized functions.

For other densities and metric structures on a plausibility simplex see
\eg\ Amari \amp\ Nagaoka \citey{amarietal1993_t2000}.

\smallskip

And measure theory? `\langitalian{La teoria della misura sta alla
  probabilit\`a come lo stucco messo male sta alle pareti: prima o poi
  cade}', said Gian-Carlo Rota
\citep{beschleretal2000,cerasoli1999,cerasolietal2003}.

\defbibnote{prenote}{See also the additional references at the end of this
  list, p.~\pageref{cha:addrefs}.

\bigskip

\small Notes:\\
Of two years separated by a virgule, the first is that of original
publication or composition.\\
 \texttt{arXiv} eprints available at \url{http://arxiv.org/},
\\
\texttt{philsci} eprints at \url{http://philsci-archive.pitt.edu/}.
}


\defbibfilter{fi}{\category{extra} \or \segment{1}}
  \printbibliography[filter=fi,prenote=prenote]
\end{refsegment}


  \section*{Additional references}
  \label{cha:addrefs}
\defbibheading{bibliography}{}

For historic and modern texts on plausibility logic, from various
perspectives, see \eg

\begin{refsegment}

  \nocite{demorgan1847,boole1854,boole1862,keynes1921_r1957,johnson1924,johnson1932,johnson1932b,johnson1932c,ramsey1926_r1950,jeffreys1931_r1957,jeffreys1939_r2003,cox1946,cox1961,hailperin1965,hailperin1976_r1986,hailperin1984,hailperin1984b,hailperin1988,hailperin1990,hailperin1996,hailperin2000,hailperin2006,hailperin2008,adams1975,adams1998,jaynes1994_r1996,pearl1988,pearl2000,portamana2007b}

\renewcommand{\citein}[2][]{\textnormal{\textcite[#1]{#2}}\addtocategory{extrb}{#2}}
\renewcommand{\citebi}[1]{\textcite{#1}\addtocategory{extrb}{#1}}
\defbibfilter{se}{\(\category{extrb} \and \not \category{extra}\) \or \segment{2}}
  \printbibliography[filter=se]
\end{refsegment}

\bigskip A glimpse on the literature presenting the \me\ principle
(excluding specific applications like statistical mechanics):

\begin{refsegment}
\nocite{jaynes1957,jaynes1963,good1963,jaynes1967,jaynes1968,bernardo1979,jaynes1979,jaynes1979b,jaynes1980c,shoreetal1980,shoreetal1981,vancampenhoutetal1981,jaynes1982,jaynes1982b,gulletal1984,csiszar1984,csiszar1985,jaynes1986b,jaynes1986d_r1996,jaynes1988b,skilling1989b,sivia1990,bucketal1991,rodriguez1991,csiszar1991,jaynes1994_r2003,aczel2004,csiszaretal2004b}

\renewcommand{\citein}[2][]{\textnormal{\textcite[#1]{#2}}\addtocategory{extrc}{#2}}
\renewcommand{\citebi}[1]{\textcite{#1}\addtocategory{extrc}{#1}}
\defbibfilter{th}{\(\category{extrc}  \and \not \category{extra} \and \not
  \category{extrb}\) \or \segment{3}}
  \printbibliography[filter=th]
\end{refsegment}
\bigskip

Discussions, criticisms, justifications and commendations of
the \me\ principle appear in

\begin{refsegment}
\nocite{hobson1969,friedmanetal1971,hobson1972,tribusetal1972,friedman1973,gageetal1973,hobsonetal1973,csiszar1975,bernardo1979b,williams1980,vanfraassen1981,csiszaretal1984,jaynes1984,jaynes1985,jaynes1985e,shimony1985,skyrms1985,jaynes1986b,seidenfeld1986_r1987,skyrms1987,zellner1988,jaynes1991c,uffink1995,uffink1996b}

\renewcommand{\citein}[2][]{\textnormal{\textcite[#1]{#2}}\addtocategory{extrd}{#2}}
\renewcommand{\citebi}[1]{\textcite{#1}\addtocategory{extrd}{#1}}
\defbibfilter{fo}{\(\category{extrd}  \and \not \category{extra} \and \not
  \category{extrb} \and \not \category{extrc} \)\or \segment{4}}
  \printbibliography[filter=fo]
\end{refsegment}

\bigskip

For some mathematical, historical, and philosophical analyses of
exchangeability and its generalizations see



\begin{refsegment}
  \nocite{johnson1924,johnson1932c,definetti1937,definetti1938,hewittetal1955,definetti1970b_t1990,heathetal1976,lindleyetal1976,diaconis1977,definetti1979,diaconisetal1980,zabell1982,jaynes1986c,diaconisetal1987,zabell1989,ladetal1990,bernardoetal1994,caves2000,caves2000c,jeffrey1980}

\renewcommand{\citein}[2][]{\textnormal{\textcite[#1]{#2}}\addtocategory{extre}{#2}}
\renewcommand{\citebi}[1]{\textcite{#1}\addtocategory{extre}{#1}}
\defbibfilter{ci}{\(\category{extre}  \and \not \category{extra} \and \not
  \category{extrb} \and \not \category{extrc} \and \not \category{extrd} \)
  \or \segment{5}}
  \printbibliography[filter=ci]

\end{refsegment}

\bigskip

Finally, for simplices and convex sets see

\begin{refsegment}
\nocite{coxeter1961_r1969,gruenbaum1967_r2003,valentine1964,broendsted1983,boydetal1004_r2009}  

\renewcommand{\citein}[2][]{\textnormal{\textcite[#1]{#2}}\addtocategory{6}{#2}}
\renewcommand{\citebi}[1]{\textcite{#1}\addtocategory{6}{#1}}
\defbibfilter{se}{\(\category{6} \and \not \category{extre}  \and \not \category{extra} \and \not
  \category{extrb} \and \not \category{extrc} \and \not \category{extrd} \)
  \or \segment{6}}
  \printbibliography[filter=se]

\end{refsegment}


\end{document}